\documentclass[12pt,notitlepage,twoside]{article}

\usepackage{latexsym}
\usepackage{amsmath}
\usepackage{amsfonts}
\usepackage{amssymb}
\let \n = \noindent
\let \dis = \displaystyle

\newcommand{\be}{\begin{equation}}
\newcommand{\ee}{\end{equation}}
\newenvironment{pf}{\noindent{\bf Proof.}\enspace}{
\hfill$\Box$\medskip}
\newenvironment{pfn}[1]{\noindent{\bf Proof of
{#1}\enspace}}{
\hfill$\Box$\medskip}

\newtheorem{thm}{Theorem}[section]
\newtheorem{prop}[thm]{Proposition}
\newtheorem{lem}[thm]{Lemma}

\newtheorem{defn}[thm]{Definition}
\numberwithin{equation}{section}

\author{Wael Abdelhedi$^1$,  Obaid Algahtani$^2$,  Hichem Chtioui$^1$ and Hichem Hajaiej$^3$ \footnote{ E-mail
 addresses: \texttt{wael\_hed@yahoo.fr} ( W. Abdelhedi),  \texttt{obalgahtani@ksu.edu.sa}  (O. AlGahtani), \texttt{Hichem.Chtioui@fss.rnu.tn} (H. Chtioui), \texttt{hichem.hajaiej@gmail.com.}(H. Hajaiej).}\\
{\footnotesize $^1$Department of mathematics,}\\ {\footnotesize
Faculty of Sciences of Sfax, 3018 Sfax,
Tunisia.}\\{\footnotesize $^2$ College of Sciences, King Saud University.}\\{\footnotesize $^3$
New York University Shanghai,
1555 Century Avenue, New Pudong, 200124 Shanghai, China.}}

\title {$Q$-curvature type problem on bounded domains of $\mathbb{R}^n$}

\begin{document}

\date{ }

\maketitle

{\footnotesize

\n{\bf Abstract.} In this paper, we establish compactness and existence results to a Branson-Paneitz type problem on a bounded domain of $\mathbb{R}^n$ with Navier boundary condition.\\
\n{\bf  MSC 2000:}\quad   35J60, 35J60, 58E05.\\
\n {\bf Key words:} Nonlinear elliptic P D E, critical exponent, Lack of compactness, Critical points at infinity.}




\section{Introduction}

In this work we are concerned with positive solutions of a nonlinear fourth order equation under the Navier boundary condition. Let $K$ be a given function on a smooth bounded domain $\Omega$ of $\mathbb{R}^n, n\geq 5$. We are looking for a map $u: \Omega\rightarrow\mathbb{R}$ satisfying the following critical fourth order PDE

\begin{equation}\label{1.1}
\left\{
  \begin{array}{c}
    \dis \Delta^2 u = \dis K(x) \;  u^\frac{n+4}{n-4},\\
    \dis u > 0 \;\mbox{ in } \Omega,\\
     \dis \Delta u= u = 0 \;\mbox{ on } \partial \Omega.
  \end{array}
\right.
\end{equation}

\n The interest of this equation comes from its resemblance to the so called $Q$-curvature problem on closed manifolds involving the Branson-Paneitz operator. The latter has widely studied in the last two decades. (See  \cite{ACR1} \cite{CR2}, \cite{BEH1}, \cite{bnc1}, \cite{AC2},  \cite{CY2}, \cite{CGY2}, \cite{DHL}, \cite{DMO3}, \cite{DMO2} and the references therein for details).

\n Problem \eqref{1.1} has a variational structure with challenging  mathematical difficulties. Indeed, if there is general standard line of attack to solve the analogous of \eqref{1.1} in the subcritical case. These  approaches do not apply to the critical case since the embedding  $\mathcal{H}\hookrightarrow L^\frac{2n}{n-4}(\Omega)$ where $\mathcal{H}:= H^2(\Omega)\cap H^1_0(\Omega)$, is not compact.

\n When $K=1$, the problem is called the Yamabe type problem. In this case, the existence of solutions of problem \eqref{1.1} depends on the topology of $\Omega$. More precisely, if $\Omega$ is a star-shaped bounded domain,  Van der Vorst \cite{V1} proved that \eqref{1.1} has no solution. When $\Omega$ has a non trivial homology group, Ebobisse-Ould Ahmedou  showed that \eqref{1.1} has  a solution \cite{EO1}.

\n When $K\neq 1$, there have been many works dealing with  \eqref{1.1}. In these contributions, the conditions on $K(x)$ ensuring the solvability of \eqref{1.1} have been discussed. In \cite{BEH1},  \cite{CB1} and \cite{CE1}, some existence results were obtained under the following two hypotheses:

$$(A) \qquad \qquad \dis \frac{\partial K}{\partial\nu}(x)\neq 0, \;\forall x\in \partial \Omega.$$
Here $\nu$ is the unit outward normal vector on $\partial \Omega$.

\n $(nd)$ \qquad$K$ is a $C^2$-positive function having only non degenerate critical points such that $$\Delta K(x)\neq 0 \mbox{ if } \nabla K(x)=0.$$

\n  Observe that $(nd)$-condition would excludes some interesting class of functions K. For example the $C^1$-functions and smooth functions having degenerate critical points. Our main motivation in this study, is to include a wider class of functions $K$ for which \eqref{1.1} admits a solution. Our main assumption is  the following $\beta$-flatness  condition:

\n $(f)_{\beta}$ \qquad Assume that $K$ is a $C^1$-positive function on $\overline{\Omega}$ such that for each critical point $y$ of $K$, there exists a real number $\beta=\beta(y)>1$, such that

$$K(x)= K(y) +\sum_{k=1}^n b_k|(x-y)_k|^\beta+ R(x-y),$$for $x$ close to $y$. Here $b_k= b_k(y)\in \mathbb{R}\setminus\{0\}$, for $k=1\ldots, n$,
and $\dis\sum_{s=0}^{[\beta]}|\nabla^sR(z)||z|^{s-\beta}=o(1)$, as $z$ tends to zero. Here $[\beta]$ denotes the integer part of $\beta$.\\

\n Note that the above mentioned $(nd)$-condition is a particular case of the $\beta$-flatness assumption (in a suitable coordinates system) taking $\beta(y)=2$ for any critical point $y$ of $K$.

\n In the first part of this paper, we are interested  with the case  $1< \beta\leq n-4, n\geq 6$. Our aim is to provide a full description of the lack of compactness of the associated variational problem to \eqref{1.1}. Indeed, we will give a characterization of all  critical points at infinity of the functional $J$ in $\Sigma^+$  and state an Euler-Hopf type of existence result.

\n Let  $G$ denote the Green's function of the bilaplacian under Navier boundary condition on $\Omega$.  It is defined by  $$G(x, y)= |x-y|^{-(n-4)}-H(x, y), \mbox{ for } x\neq y\in \Omega,$$where $H$ its regular part.\\

\n Let $\mathcal{K}$ denote the set consisting of all critical points of $K(x)$. For any $y\in \mathcal{K}$, we  define

$$\widetilde{i}(y)= \sharp\{b_k(y), \; b_k(y)<0\}.$$Let

$$\mathcal{K}_{< n-4} =\{y\in \mathcal{K}, \;\beta(y)<n-4\}, \qquad  \mathcal{K}_{< n-4}^+ =\{y\in \mathcal{K}_{< n-4}, -\sum_{k=1}^nb_k(y)>0\},$$

$$\mathcal{K}_{n-4} =\{y\in \mathcal{K}, \;\beta(y)=n-4\},  \qquad  \mathcal{K}_{n-4}^+ =\{y\in \mathcal{K}_{n-4}, -c_1\sum_{k=1}^nb_k(y)+c_2 H(y, y)>0\},$$
and $$\mathcal{K}_{>n-4} =\{y\in \mathcal{K}, \;\beta(y) >n-4\},$$
where

$$c_1= \dis \int_{\mathbb{R}^n}\frac{|x_1|^{n-4}}{(1+|x|^2)^{n}} dx,\qquad  c_2= \dis \int_{\mathbb{R}^n}\frac{dx}{(1+|z|^2)^\frac{n+4}{2}}. $$Here  $x_1$ is the first component of $x$ in some geodesic normal coordinates system.\\

\n To any $p$-tuple of distinct points $\tau_p= (y_{\ell_1}, \ldots, y_{\ell_p})\in \mathcal{K}_{n-4}^p, 1\leq p\leq \sharp \mathcal{K}$, we associate a $p\times p$ symmetric matrix $M(\tau_p)= (m_{ij})_{1\leq i, j\leq p}$ defined by:

$$m_{ii}= -\frac{1}{K(y_{\ell_i})^\frac{n}{4}}\Big(c_1\sum_{k=1}^nb_k(y_{\ell_i}) -c_2 H(y_{\ell_i}, y_{\ell_i})\Big) \mbox{ and }$$

$$m_{ij}= -c_2 \frac{G(y_{\ell_i}, y_{\ell_j})}{\Big(K(y_{\ell_i})K(y_{\ell_j})\Big)^\frac{n-4}{8}}, \mbox{ for } i\neq j.$$
Let $\rho(\tau_p)$ be the least eigenvalue of  $M(\tau_p)$.

\n $(B)$ \qquad Assume that $\dis \rho(\tau_p)\neq 0$ for any $1\leq p\leq \sharp\mathcal{K}$.

\n Lastly define

$$\mathcal C_{n-4}^{\infty}:=\big\{\tau_{p}=(y_{l_{1}},...,y_{l_{p}})\in (\mathcal K_{n-4}^+)^{p}, 1\leq p\leq\sharp\mathcal{K}, s.t.\hskip0.2cm y_{\ell_i}\neq y_{\ell_j }\hskip0.2cm \forall i\neq j, \mbox{ and }
 \rho(\tau_{p})>0\big\},$$
$$\mathcal C_{<n-4}^{\infty}:=\big\{\tau_{p}=(y_{l_{1}},...,y_{l_{p}})\in (\mathcal K_{<n-4}^+)^{p}, 1\leq p\leq\sharp\mathcal{K}, s.t.\hskip0.2cm y_{\ell_i}\neq y_{\ell_j }\hskip0.2cm \forall i\neq j\big\}$$and $$\mathcal C^\infty:= \mathcal C_{n-4}^{\infty}\cup \mathcal C_{<n-4}^{\infty}\cup \Big(\mathcal C_{n-4}^{\infty}\times\mathcal C_{<n-4}^{\infty}\Big).$$

\n The following result describes the lack of compactness of the problem \eqref{1.1}.

\begin{thm}\label{th1}
Under the assumptions   $(A), (B)$ and $(f)_{\beta}$ for $1<\beta\leq n-4$. The critical points at infinity of the associated variational problem to \eqref{1.1} ( see Definition \ref{def1}) are:

$$(y_{\ell_1}, \ldots, y_{\ell_p})_\infty:= \sum_{i=1}^{p}  \frac{1}{K(y_{\ell_i})^\frac{n-4}{2}}  P \delta_{(y_{l_{i}}, \infty)}$$
where $(y_{l_{1}}, \ldots ,y_{l_{p}})= (\tau_p) \in \mathcal{C}^\infty$. The index of a such critical points at infinity is $i(\tau_p)= p-1-\sum_{i=1}^p(n-\widetilde{i}(y_{\ell_i}))$.
\end{thm}

\n The characterization of the critical points at infinity allows us to prove the following existence result.

\begin{thm}\label{th2}
Suppose  that  $(A), (B)$ and $(f)_{\beta}$ for $1<\beta\leq n-4$ hold. If additionally

$$\sum_{\tau_p\in \mathcal C^\infty }(-1)^{i(\tau_p)}\neq 1,$$then \eqref{1.1} has a solution.
\end{thm}

\n In the second part of this paper,  we are interested to the case  of any   $\beta>1$. We prove a partial description of the lack of compactness of the problem in that case and we provide a perturbation result.

\begin{thm}\label{th3}
Assume  that $K$ satisfies $(A)$ and $(f)_{\beta}$ for $\beta>1$.  The critical points at infinity in $V(1, \varepsilon)$ are
$$ P \delta_{(y, \infty)}, \;  y\in \mathcal{K}^+_{< n-4}\cup \mathcal{K}^+_{n-4}\cup \mathcal{K}_{>n-4}.$$
Such critical points at infinity has an index  equal to $n-\widetilde{i}(y)$
\end{thm}
Now we state our perturbation result.

\begin{thm}\label{th4}
Under the assumptions  $(A)$ and $(f)_{\beta}$ for $\beta>1$, if
 $$\sum_{\dis y\in \mathcal{K}_{< n-4}^+ \; \cup \;\mathcal{K}^+_{n-4}\cup \;\mathcal{K}_{>n-4}}(-1)^{n-\widetilde{i}(y)}\neq \chi(\Omega)$$then \eqref{1.1} has a solution. Here $\chi(\Omega)$ is the Euler-Poincar\'{e} characteristic of $\Omega$.

\end{thm}

\n Our method hinges on  the  critical points at infinity theory of  A. Bahri \cite{b1}. In section 2, we state the variational structure associated to problem \eqref{1.1}. In section 3, we provide an asymptotic expansion of the gradient of $J$, without assuming any upper bound condition on the $\beta$-flatness condition. In section 4, we characterize the critical points at infinity and we prove Theorems \ref{th1} and \ref{th3}. Lastly in section 5, we prove Theorems \ref{th2} and \ref{th4}.\\

\section{Preliminaries tools}

Let $\mathcal{H}:= H^2(\Omega)\cap H^1_0(\Omega)$ with the norm $$\|u\|_{\mathcal{H}}= \Big(\int_{\Omega}(\Delta u(x))^2 dx\Big)^\frac{1}{2}.$$Define

$$ \Sigma=\bigl\{\,u\in \mathcal{H}\,,\
\mathrm{s.t.}\  \|u\|_{\mathcal{H}}=1\
  \bigr\}, \mbox{ and } \Sigma^+=\big\{\,u\in\sum\,,\
  u>0\,\bigr\}.$$
Let

$$
 J(u)=\frac{\Big(\int_{\Omega}  (\Delta u)^2\Big)^\frac{n}{n-4}}
    { \int_{\Omega}K(x)\;  u^{\frac{2n}{n-4}}dx }.
 $$

\n Observe that if $u$ is a critical point of $J$ on
$ \Sigma^+$, then $J(u)^\frac{n-4}{8}.u$ is a solution of \eqref{1.1}.\\

\n  $J$ does not satisfy the Palais-Smale
condition on $ \Sigma^+$ (P.S for short). This is due to the loss of compactness of the embedding $\mathcal{H}\hookrightarrow L^\frac{2n}{n-4}(\Omega)$. Next, we describe the sequences failing P.S condition.  For $a\in \Omega$ and $\lambda>0$, let

\begin{equation}\label{2.1}
 \delta_{a,\lambda}(x)=c_n
  \Bigl(\frac{\lambda}{1+\lambda^2|x-a|^2}\Bigr)^\frac{n-4}{2},
 \end{equation}
 where $c_n$ is a positive constant chosen such that
 $\delta_{a,\lambda}$ is the family of solutions of the following
 problem (see \cite{LI}):
  \begin{equation}
   \Delta^2u=|u|^\frac{8}{n-4}u\ ,\quad u>0\quad\mathrm{in}\
   \mathbb{R}^n.
    \end{equation}

\n Let  $P\delta_{a,\lambda}$ the unique solution of
  $$
   \left\lbrace
\begin{aligned}
\Delta^2 P\delta_{a,\lambda}&= {\delta_{a,\lambda}}^\frac{n+4}{n-4} \qquad\quad\quad     \mbox{in} \   \Omega\\
 P\delta_{a,\lambda}&=\Delta P\delta_{a,\lambda}= 0\qquad \      \mbox{on} \   \partial\Omega.
\end{aligned}
\right.
  $$

\n We have the following estimates where originally introduced by Bahri \cite{b1}.

\begin{equation}\label{2.2}
    P \delta_{a, \lambda}= \dis \delta_{a, \lambda} -\frac{c}{\lambda^{\frac{n-4}{2}}} H(a, .) + O\Big(\frac{1}{\lambda^{\frac{n}{2}} d^{n-2}}\Big),
\end{equation}
\begin{equation}\label{2.3}
   \lambda \frac{\partial P \delta_{a, \lambda}}{\partial\lambda}= \dis \lambda \frac{\partial \delta_{a, \lambda}}{\partial\lambda} + \frac{n-4}{2}\frac{c}{\lambda^{\frac{n-4}{2}}} H(a, .) + O\Big(\frac{1}{\lambda^{\frac{n}{2}} d^{n-2}}\Big),
\end{equation}
\begin{equation}\label{2.4}
 \dis \frac{1}{\lambda}\frac{\partial P \delta_{a, \lambda}}{\partial a_k}= \dis \frac{1}{\lambda} \frac{\partial \delta_{a, \lambda}}{\partial a_k} - \frac{c}{\lambda^{\frac{n-4}{2}}}\frac{\partial H(a, .)}{\partial a_k} + O\Big(\frac{1}{\lambda^\frac{n+2}{2} d^{n-1}}\Big),
\end{equation}
where $c$ is a fixed positive constant, $d= d(a, \partial\Omega)$, $a_k$ is the $k^{th}$ coordinate of $a$.\\

\n We define now the set of potential critical points at infinity associated to   $J_{\varepsilon_0}$.   Let for $\varepsilon>0$ and
$p \in \mathbb{N}^*$,
\\$$V(p,\varepsilon)=\begin{cases}u \in\Sigma^+, s.t, \hskip 0.2 cm\exists \hskip 0.2 cm a_{1},...,a_{p}\in
\Omega, \exists \;\displaystyle\lambda_{1},...,\displaystyle\lambda_{p}>\varepsilon^{-1}\; \mbox{ and }\\
\displaystyle\alpha_{1},...,\displaystyle\alpha_{p}>0
\hskip 0.2 cm with
\parallel u- \sum_{i=1}^{p}\displaystyle\alpha_{i}\displaystyle P\delta_{\displaystyle
a_{i},\displaystyle\lambda_{i}}\parallel<\varepsilon, \; \varepsilon_{ij}<\varepsilon \hskip 0.2 cm\forall i\neq j,\\
\lambda_i d_i > \varepsilon^{-1}\; \mbox{ and } \big|J^{\frac{n}{n-4}}(u)\displaystyle\alpha_{i}^\frac{8}{n-4}K(a_{i})-1\big|<\varepsilon
\hskip 0.2 cm \forall i=1, \ldots, p.
\end{cases}$$Here, $d_i= d(a_i, \partial \Omega)$ and  $\varepsilon_{ij}=\biggr(\displaystyle
\frac{\displaystyle\lambda_{i}}{\displaystyle\lambda_{j}}+\displaystyle
\frac{\displaystyle\lambda_{j}}{\displaystyle\lambda_{i}}+\displaystyle
\lambda_{i}\lambda_{j} |a_i-a_j|^2 \biggr)^{\frac{4-n}{2}}.$\\

\begin{prop}\label{p2.1}(\cite{bc2}, \cite{S3}) Assume that $J$ has no critical points in $\Sigma^+$.  Let $(u_k)_k$ be a sequence in $\Sigma^+$ such that
$J (u_k)$ is bounded and $\partial J (u_k)$ goes to zero. Then there
exists a positive  integer $p$, a sequence $(\varepsilon_k)
$ with  $\varepsilon_k\rightarrow 0$ as $k\rightarrow +\infty$  and an extracted subsequence
of $(u_k)_k$'s, again denoted $(u_k)_k$, such that $u_k \in
V(p,\varepsilon_k), \forall k$.
\end{prop}

\n The following Proposition gives a parametrization of $V(p, \varepsilon)$.

\begin{prop}\label{p2.2}(\cite{bc2})
 For all  $p \in \mathbb{N}^*$, there exists  $\varepsilon_{p}>0$ such that for any
$\varepsilon\leq\varepsilon_{p}$ and any $u$ in $ V(p,\varepsilon)$,
the problem

$$\min \Big\{\Big\| u
-\dis\sum_{i=1}^{p}
\alpha_{i} P \delta_{a_{i},\lambda_{i}}\Big\|, \;\alpha_{i}>0, \lambda_{i}>0, a_{i}\in
\Omega \Big\}.$$
has a unique solution  (up to a permutation). Thus, we can uniquely write $u$ as follows

$$u=\dis\sum_{i=1}^{p}\alpha_{i} P \delta_{a_{i},
\lambda_{i}}+v,$$ where $v\in
H^2(\Omega)\cap H^1_0(\Omega)\cap T_wW_s(w)$  and satisfies
$$
(V_{0}): \Big\langle v ,\psi \Big\rangle=0 \; for\;
    \psi\in\{\displaystyle P\delta_{i}, \frac{\displaystyle\partial P\delta_{i}}{\displaystyle\partial\lambda_{i}}, \frac{\displaystyle\partial P\delta_{i}}
    {\displaystyle\partial
    a_{i}},i=1,...,p\}.$$Here, $P\delta_i
= P\delta_{a_{i},\lambda_{i}}$ and $<.,.>$ denotes the inner
product  on  $H^2(\Omega)\cap H^1_0(\Omega)$ defined by
  $$\langle u \,,\,v  \rangle =\int_\Omega\Delta u \Delta v.$$
\end{prop}

\n The following Proposition deals with the $v$-part of $u$ and shows that is negligible with respect to  the concentration phenomenon.

\begin{prop} \label{}(\cite{bc2}, \cite{b1})
There is a $\mathcal{C}^{1}$-map which to each
$(\alpha_{i},a_{i},\lambda_{i})$ such that\\ $\dis\sum_{i=1}^{p}
\alpha_{i} P\delta_{a_{i},\lambda_{i}}$ belongs to
$V(p, \varepsilon)$ associates $\overline{v}=\overline
v(\alpha_i, a_i, \lambda_i)$ such that $\overline{v}$ is the unique solution of the following minimization problem

$$\dis\min \Big\{J\Big(\sum_{i=1}^{p}
\alpha_{i} P\delta_{a_{i},\lambda_{i}}+v\Big), v \in H^1_0(\Omega)\; \mbox{ and satisfies } (V_0)\Big\}.$$
Moreover, there exists a change of variables
$v-\overline{v}\rightarrow V$ such that

$$J \Big(\dis\sum_{i=1}^{p}
\alpha_{i} P\delta_{a_{i},\lambda_{i}}+v\Big) = J \Big(\dis\sum_{i=1}^{p}
\alpha_{i} P\delta_{a_{i},\lambda_{i}}+\overline{v}\Big)+\parallel
V \parallel^{2}.$$
\end{prop}

\n We now state the definition of critical point at infinity.

\begin{defn}\label{def1}\cite{b1}
A critical point at infinity of $J$ is a limit of a non-compact
flow line $u(s)$ of the gradient vector field $(-\partial J)$. By Propositions \ref{p2.1} and \ref{p2.2}, $u(s)$ can be
written as:
\begin{center}
$u(s)=\dis \sum_{i=1}^{p}\dis\alpha_{i}(s)\dis P\delta_{\dis
a_{i}(s),\dis\lambda_{i}(s)}+v(s)$.
\end{center}
Denoting by ${y}_{i} = \dis\lim_{s \longrightarrow + \infty} \dis
a_{i}(s)\;\; \mbox{ and }\;  \dis{\alpha}_{i} =\dis \lim_{s \longrightarrow
+\infty} \dis\alpha_{i}(s)
$, we then denote by

\begin{center}
$ \dis
\sum_{i=1}^{p}\dis \alpha_{i}\dis P\delta_{
y_{i},\infty}$ or $(\dis
y_{1},...,\dis y_{p})_{\infty}$
\end{center} such a critical point at
infinity.
\end{defn}

\section{Expansion of the gradient of $J$}

Let $\rho$ be a positive small constant such that for any $y\in \mathcal{K}$, the expansion $(f)_\beta$ holds in $B(y, \rho)$. Let
$$\widetilde{V}(p, \varepsilon):= \{u= \sum_{i=1}^p\alpha_i P \delta_{(a_i, \lambda_i)}\in {V}(p, \varepsilon), \; a_i\in B(y_{\ell_i}, \rho), y_{\ell_i}\in \mathcal{K}, \forall i=1, \ldots, p\}.$$
The following proposition gives the variation of $J$ in $\widetilde{V}(p, \varepsilon)$ with respect to $\lambda_i, i=1, \ldots, p$.

\begin{prop}\label{lem1}
Assume  that $K$ satisfies  $(f)_{\beta}$ condition for $\beta>1$. For  $u=\sum_{i=1}^p  \alpha_i P \delta_{a_i, \lambda_i}\in \widetilde{V}(p, \varepsilon)$, we have the following two estimates:

\begin{eqnarray}
 \nonumber
&(a)& \; \Big \langle \partial J(u), \alpha_i\lambda_{i} \dis \frac{\partial
P \delta_{a_i, \lambda_i}}{\partial \lambda_{i}} \Big\rangle =
 - 2c_2 \frac{J(u)}{K(a_i)}
    \dis \sum_{j \neq i} \alpha_i\alpha_{j} \Big(\lambda_{i}\dis \frac{\partial \varepsilon_{i j}}{\partial
    \lambda_{i}} + \frac{H(a_i, a_j)}{(\lambda_i\lambda_j)^\frac{n-4}{2}}\Big)\\ \nonumber &+&  O\Big(\sum_{j=2}^{[\beta]}\frac{|a_i- y_{\ell_i}|^{\beta-j}}{\lambda_i^j}\Big)
      + \dis O\Big( \frac{1}{\lambda_{i}^{\beta}}\Big) + o\bigg( \dis \sum_{j \neq i} \Big(\varepsilon_{i j}+ \frac{H(a_i, a_j)}{(\lambda_i\lambda_j)^\frac{n-4}{2}}
    \Big)\bigg).\\
 \nonumber
&(b)& \; \Big \langle \partial J(u), \alpha_i\lambda_{i} \dis \frac{\partial
P \delta_{a_i, \lambda_i}}{\partial \lambda_{i}} \Big\rangle =
 - 2c_2 \frac{J(u)}{K(a_i)}
    \dis \sum_{j \neq i} \alpha_i\alpha_{j} \Big(\lambda_{i}\dis \frac{\partial \varepsilon_{i j}}{\partial
    \lambda_{i}} + \frac{H(a_i, a_j)}{(\lambda_i\lambda_j)^\frac{n-4}{2}}\Big)\\ \nonumber &+& 2\alpha_i^2\frac{J(u)}{K(a_i)}
    \left\{
      \begin{array}{ll}
           \dis \frac{n-4}{2}c_1\frac{\sum_{k=1}^n b_k(y_{\ell_i})}{\lambda_i^{\beta(y_{\ell_i})}}, &\mbox{ if } \beta(y_{\ell_i})< n-4 \\
            \dis \frac{n-4}{2}c_1\frac{\sum_{k=1}^n b_k(y_{\ell_i})}{\lambda_i^{\beta(y_{\ell_i})}}-c_2\frac{H(y_{\ell_i}, y_{\ell_i})}{\lambda_i^{n-4}}, &\mbox{ if }  \beta(y_{\ell_i})= n-4 \\
\dis -c_2\frac{H(y_{\ell_i}, y_{\ell_i})}{\lambda_i^{n-4}}, &\mbox{ if }  \beta(y_{\ell_i})>n-4
      \end{array}
    \right.
     \\  \nonumber &+& O\Big(|a_i- y_{\ell_i}|^{\beta}\Big) +o\bigg( \dis \sum_{j \neq i} \Big(\varepsilon_{i j}+ \frac{H(a_i, a_j)}{(\lambda_i\lambda_j)^\frac{n-4}{2}}
    \Big) + \frac{1}{\lambda_{i}^{n-4}}\bigg).
\end{eqnarray}
Here  $c_1= \dis \int_{\mathbb{R}^n}|z_1|^{\beta}\frac{|z|^2-1}{(1+|z|^2)^{n+1}} dz$ and $c_2= \dis \int_{\mathbb{R}^n}\frac{|z|^2-1}{(1+|z|^2)^{n}} dz$.

\end{prop}

\begin{pf}
For  $u=\sum_{i=1}^p  \alpha_i P \delta_{a_i, \lambda_i}\in \widetilde{V}(p, \varepsilon)$, we have:

$$\partial J(u)= 2J(u) \Big(u+ J(u) \Delta^{-1}(K u^\frac{n+4}{n-4})\Big).$$Thus,

$$\Big \langle \partial J(u),  \alpha_i\lambda_{i} \dis \frac{\partial
P \delta_{a_i, \lambda_i}}{\partial \lambda_{i}} \Big \rangle = 2 J(u) \bigg[\sum_{j=1}^p \alpha_i\alpha_j<
P \delta_{a_j, \lambda_j},  \lambda_{i} \dis \frac{\partial
P \delta_{a_i, \lambda_i}}{\partial \lambda_{i}} >$$$$ -J(u)^\frac{n}{n-4} \int_{\Omega} K (\sum_{j=1}^p \alpha_j P \delta_{a_j, \lambda_j})^\frac{n+4}{n-4} \alpha_i\lambda_i  \dis \frac{\partial
P \delta_{a_i, \lambda_i}}{\partial \lambda_{i}} \bigg].$$

\n Using \eqref{2.2} and \eqref{2.3} and the fact that $\dis J(u)^\frac{n}{n-4}\alpha_i^\frac{4}{n-4}K(a_i)=1+o(1)$, we get

$$\Big \langle \partial J(u),  \alpha_i\lambda_{i} \dis \frac{\partial
P \delta_{a_i, \lambda_i}}{\partial \lambda_{i}} \Big \rangle = - 2c_2 \frac{J(u)}{K(a_i)}
    \dis \sum_{j \neq i} \alpha_i\alpha_{j} \Big(\lambda_{i}\dis \frac{\partial \varepsilon_{i j}}{\partial
    \lambda_{i}} + \frac{H(a_i, a_j)}{(\lambda_i\lambda_j)^\frac{n-4}{2}}\Big)$$$$ -2\alpha_i^2\frac{J(u)}{K(a_i)}\bigg(\int_{\Omega} K \delta_{a_j, \lambda_j}^\frac{n+4}{n-4} \lambda_i  \dis \frac{\partial
 \delta_{a_i, \lambda_i}}{\partial \lambda_{i}}+ c_2\frac{H(a_i, a_i)}{\lambda_i^{n-4}} \bigg) + o\bigg( \dis \sum_{j \neq i} \Big(\varepsilon_{i j}+ \frac{H(a_i, a_j)}{(\lambda_i\lambda_j)^\frac{n-4}{2}}
    \Big)+ \frac{1}{\lambda_i^{n-4}}\bigg).$$Observe that

$$\delta_{a, \lambda}^\frac{n+4}{n-4} \lambda \frac{\partial \delta_{a, \lambda}}{\partial \lambda}= \frac{n-4}{2} \lambda^n\dis\frac{1-\lambda^2|x-a|^2}{(1+\lambda^2|x-a|^2)^{n+1}}.$$

\n Let $\mu>0$  such that $B(a_i, \mu)\subset B(y_{\ell_i}, \rho)$. We have

$$ \int_{\Omega} K(x) \delta_{a_i, \lambda_i}^\frac{n+4}{n-4} \lambda \frac{\partial \delta_{a_i, \lambda_i}}{\partial \lambda_i} dx=  \int_{B(a_i, \mu)} K(x) \delta_{a_i, \lambda_i}^\frac{n+4}{n-4} \lambda_i \frac{\partial \delta_{a_i, \lambda_i}}{\partial \lambda_i} dx +O\Big(\frac{1}{\lambda_i^n}\Big).$$
 After a change of variables $z= \lambda_i(x-a_i)$,

$$ \int_{B(a_i, \mu)} K(x) \delta_{a_i, \lambda_i}^\frac{n+4}{n-4} \lambda_i \frac{\partial \delta_{a_i, \lambda_i}}{\partial \lambda_i} dx=\frac{n-4}{2} \int_{B(0, \lambda_i\mu)} K(a_i+ \frac{z}{\lambda_i})  \dis\frac{1-|z|^2}{(1+|z|^2)^{n+1}}dz.$$
Using the following expansion of $K$  around $a_i$,

$$ K(a_i+ \frac{z}{\lambda_i}) = K(a_i)+ \sum_{j=1}^{[\beta]}\frac{D^j K(a_i)(\frac{z}{\lambda_i})^j}{j!}+ O\Big(|\frac{z}{\lambda_i}|^\beta\Big),$$
and the fact that $\dis \int_{B(0, \lambda_i\mu)} K(a_i)  \dis\frac{1-|z|^2}{(1+|z|^2)^{n+1}}dz= O\Big(\frac{1}{\lambda_i^n}\Big)$, we get

$$ \int_{B(a_i, \mu)} K(x) \delta_{a_i, \lambda_i}^\frac{n+4}{n-4} \lambda_i \frac{\partial \delta_{a_i, \lambda_i}}{\partial \lambda_i} dx=
\frac{n-4}{2} \sum_{j=1}^{[\beta]}  \int_{B(0, \lambda_i\mu)}  \dis \frac{D^j K(a_i)(z)^j}{\lambda_i^jj!} \frac{1-|z|^2}{(1+|z|^2)^{n+1}}dz$$$$ + O\bigg(\frac{1}{\lambda_i^\beta}\int_{\mathbb{R}^n}  \dis |z_k|^\beta\frac{1-|z|^2}{(1+|z|^2)^{n+1}}dz\bigg) +O\bigg(\frac{1}{\lambda_i^{n}}\bigg).$$

\n Observe that,
$$\int_{B(0, \lambda_i\mu)}  \dis \dis {D K(a_i)(z)} \frac{1-|z|^2}{(1+|z|^2)^{n+1}} dz = 0.$$
Moreover, under $(f)_{\beta}$-condition, we have

$$|D^jK(a_i)|= O(|a_i- y_{\ell_i}|^{\beta-j}).$$Thus,

$$\int_{B(a_i, \mu)} K(x) \delta_{a_i, \lambda_i}^\frac{n+4}{n-4} \lambda_i \frac{\partial \delta_{a_i, \lambda_i}}{\partial \lambda_i} dx=
 O\bigg( \sum_{j=2}^{[\beta]}   \dis \frac{|a_i- y_{\ell_i}|^{\beta-j})}{\lambda_i^j}\bigg) +O\bigg(\frac{1}{\lambda_i^{\beta}}\bigg).$$Hence, the estimate $(a)$ of Proposition  \ref{lem1} follows.

\n For the estimate $(b)$,  $(f)_{\beta}$-expansion yields

$$K(a_i + \frac{z}{\lambda_i})= K(y_{\ell_i}) + \sum_{k=1}^n b_k \Big|  \frac{z_k}{\lambda_i}+ (a_i-y_{\ell_i})_k\Big|^\beta + o\bigg(\bigg| \frac{z}{\lambda_i}+ (a_i-y_{\ell_i})\bigg|^\beta\bigg)$$Therefore,

$$\int_{B(a_i, \mu)} K(x) \delta_{a_i, \lambda_i}^\frac{n+4}{n-4} \lambda_i \frac{\partial \delta_{a_i, \lambda_i}}{\partial \lambda_i} dx= \dis\frac{n-4}{2}\frac{1}{\lambda_i^\beta} \sum_{k=1}^n b_k  \int_{B(0, \lambda_i\mu)}  \Big|z_k+ {\lambda_i}(a_i-y_{\ell_i})_k\Big|^\beta \dis\frac{1-|z|^2}{(1+|z|^2)^{n+1}}dz$$$$ + o\bigg(\frac{1}{\lambda_i^\beta} \int_{\mathbb{R}^n}  |z_k |^\beta \dis\frac{1-|z|^2}{(1+|z|^2)^{n+1}}dz\bigg)+
o\bigg( |a_i-y_{\ell_i}|^\beta \int_{\mathbb{R}^n}  |z_k |^\beta \dis\frac{1-|z|^2}{(1+|z|^2)^{n+1}}dz\bigg).$$
$$= \dis\frac{n-4}{2}\frac{1}{\lambda_i^\beta} \sum_{k=1}^n b_k  \int_{B(0, \lambda_i\mu)}  \Big|z_k\Big|^\beta \dis\frac{1-|z|^2}{(1+|z|^2)^{n+1}}dz+ O( |a_i-y_{\ell_i}|^\beta)$$$$ + o\bigg(\frac{1}{\lambda_i^\beta} \int_{B(0, \lambda_i\mu)}  |z_k |^\beta \dis\frac{|1-|z|^2|}{(1+|z|^2)^{n+1}}dz\bigg)+ o\bigg(\frac{1}{\lambda_i^{n-4}}\bigg).$$

\n Observe that, for $\beta < n$,

$$\int_{B(0, \lambda_i\mu)}  \dis |z_k|^\beta \frac{1-|z|^2}{(1+|z|^2)^{n+1}} dz =  \int_{\mathbb{R}^n}  \dis |z_k|^\beta \frac{1-|z|^2}{(1+|z|^2)^{n+1}} dz + O(\frac{1}{\lambda_i^{n-\beta}})=-c_1+ O(\frac{1}{\lambda_i^{n-\beta}}).$$
For  $\beta=n$,

$$ \int_{B(0, \lambda_i\mu)}  \dis |z_k|^\beta \frac{1-|z|^2}{(1+|z|^2)^{n+1}} dz =O(\log \lambda_i).$$

\n Lastly,  for $\beta>n$,   $$\int_{B(0, \lambda_i\mu)}  \dis |z_k|^\beta \frac{1-|z|^2}{(1+|z|^2)^{n+1}} dz =O(\frac{1}{\lambda_i^{n-\beta}}).$$Therefore,

\begin{equation}\label{}
    \int_{\Omega} K(x) \delta_{a, \lambda}^\frac{n+4}{n-4} \lambda \frac{\partial \delta_{a, \lambda}}{\partial \lambda} dx= \bigg(-\dis \frac{n-4}{2} c_1 \frac{\sum_{k=1}^n b_k}{\lambda_i^\beta} + o(\frac{1}{\lambda_i^{\beta}}), \mbox{ if } \beta< n\bigg)+ O( |a_i-y_{\ell_i}|^\beta) + o (\frac{1}{\lambda_i^{n-4}}).
\end{equation}
This conclude   the proof of Proposition \ref{lem1}.
\end{pf}

\begin{prop}\label{lem2}
Assume  that $K$ satisfies  $(f)_{\beta}$ condition for $\beta>1$. Let  $u=\dis\sum_{i=1}^p  \alpha_i P \delta_{a_i, \lambda_i}\in \widetilde{V}(p, \varepsilon)$. For any $i=1, \ldots, p$ and $k=1, \ldots, n$, we have the  following expansions.

\begin{eqnarray}
 \nonumber
(i) &&\Big \langle \partial J(u), \alpha_i\frac{1}{\lambda_{i}} \dis \frac{\partial
P\delta_{a_i, \lambda_i}}{\partial (a_{i})_k}  \Big\rangle =
  - c_3 \alpha_i^2J(u) \frac{b_k}{\lambda_i K(a_i)} \beta\mbox{ sign}(a_i-y_{\ell_i})_k |(a_i-y_{\ell_i})_k|^{\beta-1}\\
     &+& O\Big(\sum_{j=2}^{\min (n, \beta)}\frac{|a_i- y_{\ell_i}|^{\beta-j}}{\lambda_i^j}\Big)
+    O\Big(\frac{1}{\lambda_i^{\min (n, \beta)}}\Big) +    O\Big(\frac{1}{\lambda_i^{n-1}}\Big)+ O\Big( \dis \sum_{j \neq i} |\frac{1}{\lambda_i}\frac{\partial \varepsilon_{i j}}{\partial a_i}|
    \Big).
\end{eqnarray}
Moreover, if $\lambda_{i} |a_{i}-y_{\ell_{i}}|$ is bounded and $\beta< n+1$, we have

\begin{eqnarray}
 \nonumber
(ii) \Big \langle  \partial J(u), \alpha_i\frac{1}{\lambda_{i}} \dis \frac{\partial
P\delta_{a_i, \lambda_i}}{\partial (a_{i})_k}  \Big\rangle &=&  - (n-2)\alpha_i^2J(u) \frac{b_k}{K(a_i)\lambda_i^\beta}\int_{\mathbb{R}^{n}}|z_k+ \lambda_i(a_i-y_{\ell_i})_k|^{\beta}\\
     &\times&\frac{z_k}{(1+|z|^2)^{n+1}}dz +   o\Big(\frac{1}{\lambda_i^{\beta}}\Big) + O\Big( \dis \sum_{j \neq i} |\frac{1}{\lambda_i}\frac{\partial \varepsilon_{i j}}{\partial a_i}|
    \Big).
\end{eqnarray}
Here $c_3= \dis ((n-2)c_n^\frac{2n}{n-4})\int_{\mathbb{R}^n}  \dis \frac{|z|^2}{(1+|z|^2)^{n+1}} dz$.

\end{prop}

\begin{pf}
We argue as in the proof of Proposition  \ref{lem1},

$$\Big \langle \partial J(u), \alpha\frac{1}{\lambda}\frac{\partial P \delta_{a, \lambda} }{\partial a_k}\Big \rangle = \frac{2\alpha^2 J(u)}{K(a)}\bigg( -\int_{\Omega} K(x) \delta_{a, \lambda}^\frac{n+4}{n-4} \frac{1}{\lambda} \frac{\partial \delta_{a, \lambda}}{\partial a_k} dx + O\bigg(\frac{\frac{\partial H}{\partial a_k}(a, a)}{\lambda^{n-1}}\bigg)  +o\bigg(\frac{1}{\lambda^{n-1}}\bigg)\bigg)$$$$+  O\Big( \dis \sum_{j \neq i} |\frac{1}{\lambda_i}\frac{\partial \varepsilon_{i j}}{\partial a_i}|\Big)$$$$= -\frac{2\alpha^2 J(u)}{K(a)} \int_{B(a, \mu)} K(x) \delta_{a, \lambda}^\frac{n+4}{n-4} \frac{1}{\lambda} \frac{\partial \delta_{a, \lambda}}{\partial a_k} dx + O\bigg(\frac{ 1}{\lambda^{n-1}}\bigg)+ O\Big( \dis \sum_{j \neq i} |\frac{1}{\lambda_i}\frac{\partial \varepsilon_{i j}}{\partial a_i}| \Big).$$Observe that

$$ \delta_{a, \lambda}^\frac{n+4}{n-4} \frac{1}{\lambda} \frac{\partial \delta_{a, \lambda}}{\partial a_k} = (n-4)c_n^\frac{2n}{n-4} \frac{\lambda^{n+1}(x-a)_k}{(1+\lambda^2|x-a|^2)^{n+1}}.$$
A change of variables $z= \lambda(x-a)$ yields

$$  \int_{B(a, \mu)} K(x) \delta_{a, \lambda}^\frac{n+4}{n-4} \frac{1}{\lambda} \frac{\partial \delta_{a, \lambda}}{\partial a_k} dx = (n-4)c_n^\frac{2n}{n-4} \int_{B(0, \lambda\mu)} K(a + z/\lambda) \frac{z_k}{(1+|z|^2)^{n+1}} dz.$$
To get the first expansion of Proposition \ref{lem2}, we expand $K$ as follows $$K(a + z/\lambda) =K(a)+ \sum_{j=1}^{\min (n, \beta)}\frac{D^jK(a)(z/\lambda)^j}{j!}+ O\bigg(\frac{|z|^{\min(\beta, n)}}{\lambda^{\min(\beta, n)}}\bigg).$$Using the fact that $\dis \int_{B(0, \lambda\mu)} K(a) \frac{z_k}{(1+|z|^2)^{n+1}} dz =0$, we get

$$  \int_{B(a, \mu)} K(x) \delta_{a, \lambda}^\frac{n+4}{n-4} \frac{1}{\lambda} \frac{\partial \delta_{a, \lambda}}{\partial a_k} dx = (n-4)c_n^\frac{2n}{n-4}
\frac{1}{\lambda} \int_{B(0, \lambda\mu)}  \frac{DK(a)(z). z_k}{(1+|z|^2)^{n+1}} dz$$$$ +O\bigg(\sum_{j=2}^{\min (n, \beta)}\frac{1}{\lambda^j} \int_{\mathbb{R}^n}  \frac{|D^jK(a)||z|^{j+1}}{(1+|z|^2)^{n+1}} dz \bigg)
+ O\bigg(\frac{1}{\lambda^{\min(n, \beta)}} \int_{\mathbb{R}^n}  \frac{|z|^{\min(n, \beta)+1}}{(1+|z|^2)^{n+1}} dz \bigg).$$Observe that

$$\int_{B(0, \lambda\mu)}  \frac{DK(a)(z). z_k}{(1+|z|^2)^{n+1}} dz= \sum_{j=1}^n\frac{\partial K(a)}{\partial x_j}\int_{B(0, \lambda\mu)}  \frac{z_j z_k}{(1+|z|^2)^{n+1}} dz$$$$= \frac{1}{n} \frac{\partial K(a)}{\partial x_k}\int_{B(0, \lambda\mu)}  \frac{|z|^2}{(1+|z|^2)^{n+1}} dz=\frac{1}{n} \frac{\partial K(a)}{\partial x_k}\int_{\mathbb{R}^n}  \frac{|z|^2}{(1+|z|^2)^{n+1}} dz + O\Big(\frac{1}{\lambda^{n}}\Big),$$since $\dis \int_{B(0, \lambda\mu)}  \frac{z_j z_k}{(1+|z|^2)^{n+1}} dz=0, \forall j\neq k$. Using now the fact that $a\in B(y, \rho)$, we derive from $(f)_\beta$-condition  that

$$\frac{\partial K}{\partial x_k}(a) = b_k \beta \mbox{ sign}(a-y)_k|(a-y)_k|^{\beta-1} + \frac{\partial R}{\partial x_k}(a-y)$$
$$=b_k \beta \mbox{ sign}(a-y)_k|(a-y)_k|^{\beta-1} +o(|x-y|^{\beta-1}).$$

\n Moreover, for every $j=2, \ldots \tilde{\beta},$

$$|D^j K(a)|= O(|a -y|^{\beta-j}).$$Thus,

$$  \int_{B(a, \mu)} K(x) \delta_{a, \lambda}^\frac{n+4}{n-4} \frac{1}{\lambda} \frac{\partial \delta_{a, \lambda}}{\partial a_k} dx =c_3\frac{ b_k}{\lambda}  \mbox{ sign}(a-y)_k {|(a-y)_k|^{\beta-1}}+ o\Big(\frac{|a-y|^{\beta-1}}{\lambda}\Big)$$$$+O\bigg(\sum_{j=2}^{\tilde{\beta}}\frac{ |a -y|^{\beta-j} }{\lambda^j}\bigg)+ O\bigg(\frac{1}{\lambda^{\min(n, \beta)}}\bigg).$$This finishes  the proof of (a) of Proposition \ref{lem2}. Concerning the estimate (b), it follows from the above arguments and the following estimate

$$  \int_{B(a, \mu)} K(x) \delta_{a, \lambda}^\frac{n+4}{n-4} \frac{1}{\lambda} \frac{\partial \delta_{a, \lambda}}{\partial a_k} dx =(n-4)c_n^\frac{2n}{n-4} \frac{b_k}{\lambda}\int_{\mathbb{R}^n}  \dis |z_k+\lambda(a-y)_k|^\beta\frac{z_k}{(1+|z|^2)^{n+1}}dz + o\bigg(\frac{1}{\lambda^{\beta}}\bigg).$$
This finishes  the proof of Proposition \ref{lem2}.

\end{pf}

\section{Lack of compactness and critical points at infinity}

In the first part of this section, we focus on $\widetilde{V}(1, \varepsilon)$; the neighborhood of critical points at infinity consisting by single masses. We study the concentration phenomenon in this set and we identify the related critical points at infinity.  Let $\rho>0$ small enough such that for any $y\in \mathcal{K}$, the expansion $(f)_\beta$ holds in $B(y, \rho)$ and let:

\begin{eqnarray}
\nonumber \widetilde{V}_1(1, \varepsilon) &=& \{u= \alpha\delta_{(a, \lambda)}\in \widetilde{V}(1, \varepsilon), \; a\in B(y, \rho), \; y\in \mathcal{K} \mbox{ with } \beta= \beta(y) < n-4\},\\
\nonumber \widetilde{V}_2(1, \varepsilon) &=& \{u= \alpha\delta_{(a, \lambda)}\in \widetilde{V}(1, \varepsilon), \; a\in B(y, \rho), \; y\in \mathcal{K} \mbox{ with } \beta= \beta(y) = n-4\},\\
\nonumber \widetilde{V}_3(1, \varepsilon) &=& \{u= \alpha\delta_{(a, \lambda)}\in \widetilde{V}(1, \varepsilon), \; a\in B(y, \rho), \; y\in \mathcal{K} \mbox{ with } \beta= \beta(y) > n-4\}
\end{eqnarray}

\n As in \cite{b1}, see also \cite{ach1}, the characterization of the critical points at infinity in $\widetilde{V}(1, \varepsilon)$ is obtained through the construction of a suitable decreasing pseudo-gradient satisfying the P.S condition as long as the concentration point $a(s)$ does not enter in a neighborhood of $y\in \mathcal{K}_{< n-4}^+ \; \cup \;\mathcal{K}^+_{n-4}\cup \;\mathcal{K}_{>n-4}$.

\n let $\delta$ be a small positive constant  and let $\theta_1, \theta_2$ and $\theta_3$ be the following three cut-off functions
\begin{eqnarray}
\nonumber \theta_1: \mathbb{R}&\longrightarrow& \mathbb{R}\\
\nonumber t&\longmapsto& \left\{
                           \begin{array}{ll}
                             1 & \hbox{ if } \dis|t|\leq\frac{\delta}{2}\\
                             0 & \hbox{ if } |t| \geq {\delta}.
                           \end{array}
                         \right.
\end{eqnarray}
\begin{eqnarray}
\nonumber \theta_2: \mathbb{R}&\longrightarrow& \mathbb{R}\\
\nonumber t&\longmapsto& \left\{
                           \begin{array}{ll}
                             1 & \hbox{ if } \dis\frac{\delta}{2} \leq |t|\leq\frac{1}{\delta}\\
                             0 & \hbox{ if } |t|\dis \in [0, \frac{\delta}{4}]\cup [\frac{2}{\delta}, +\infty[.
                           \end{array}
                         \right.
\end{eqnarray}
\begin{eqnarray}
\nonumber \theta_3: \mathbb{R}&\longrightarrow& \mathbb{R}\\
\nonumber t&\longmapsto& \left\{
                           \begin{array}{ll}
                             1 & \hbox{ if } |t|\geq \dis \frac{1}{\delta}\\
                             0 & \hbox{ if } |t| \leq \dis\frac{1}{2\delta}.
                           \end{array}
                         \right.
\end{eqnarray}

\smallskip

\n\textbf{$\bullet$ Pseudo-gradient in $\widetilde{V}_1(1, \varepsilon)$:}

\n Let $W_1$ be the following vector field. $\forall u= \alpha\delta_{(a, \lambda)}\in \widetilde{ V}_1(1, \varepsilon)$,
$$W_1(u)= - \theta_1(\lambda |a-y|) (\sum_{k=1}^n b_k) \dis \alpha\lambda \frac{ \partial \delta_{(a, \lambda)}}{\partial \lambda} +\theta_3(\lambda |a-y|) (\sum_{k=1}^n b_k) \mbox{ sign } (a-y)_k \; \dis\alpha\frac{1}{\lambda}\frac{\partial \delta_{(a, \lambda)}}{\partial a_k}$$$$ +  \theta_2(\lambda |a-y|)(\sum_{k=1}^n b_k)\int_{\mathbb{R}^n}  |z_k+ \lambda(a-y)_k |^\beta\frac{ z_k}{(1 +|z|^2)^{n+1}}  dz \; \dis\alpha\frac{1}{\lambda}\frac{\partial \delta_{(a, \lambda)}}{\partial a_k}.$$
We claim that

\begin{equation}\label{aa}
\big<\partial J(u), W_1(u)\big> \leq -c \Big(\frac{1}{\lambda^{\beta}}+ \frac{|\nabla K(a)|}{\lambda}\Big).
\end{equation}
Indeed,  if $\lambda|a-y|\leq\delta$, by Proposition \ref{lem1}, we have

\begin{equation}\label{3a}
<\partial J(u), -(\sum_{k=1}^nb_k) \dis \alpha\lambda \frac{ \partial \delta_{(a, \lambda)}}{\partial \lambda}>\leq  - c \; \frac{ (\sum_{k=1}^nb_k)^2}{\lambda^\beta},
\end{equation}since $\dis|a-y|^\beta=o\big(\frac{1}{\lambda^{\beta}}\big)$ as $\delta$ small enough. Observe that under $(f_\beta)$-condition, we have

\begin{equation}\label{3a1}
    |\nabla K(a)|= O\Big(|a-y|^{\beta-1}\Big),
\end{equation}thus

\begin{equation}\label{3b}
   \frac{ |\nabla K(a)|}{\lambda}= O\Big(\frac{1}{\lambda^{\beta}}\Big), \; \; \mbox{ if }\lambda|a-y| \mbox{ is bounded}.
\end{equation}
Therefore, we can appear $\dis -  \frac{ |\nabla K(a)|}{\lambda}$ in the upper bound of \eqref{3a} and we obtain

\begin{equation}\label{}
<\partial J(u), -(\sum_{k=1}^nb_k) \dis \alpha\lambda \frac{ \partial \delta_{(a, \lambda)}}{\partial \lambda}>\leq  - c\; \bigg( \frac{1}{\lambda^\beta} + \frac{ |\nabla K(a)|}{\lambda}\bigg).
\end{equation}

\n If $\lambda|a-y| \in [\frac{\delta}{4}, \frac{2}{\delta}]$, by the second expansion of Proposition \ref{lem2}, we obtain

$$\big<\partial J(u), \sum_{k=1}^n b_k \int_{\mathbb{R}^n}  |z_k+ \lambda(a-y)_k |^\beta\frac{ z_k}{(1 +|z|^2)^{n+1}}  dz \; \dis\alpha\frac{1}{\lambda}\frac{\partial \delta_{(a, \lambda)}}{\partial a_k}\big>$$
$$ \leq -\frac{c}{\lambda^\beta} \bigg(\int_{\mathbb{R}^n}  |z_{k_a}+ \lambda(a-y)_{k_a} |^\beta\frac{ z_{k_a}}{(1 +|z|^2)^{n+1}}  dz\bigg)^2+ o(\frac{1}{\lambda^\beta}),$$
where $|(a-y)_{k_a} |= \dis\max_{1\leq k \leq n} |(a-y)_{k} |$. Since $\lambda|(a-y)_{k_a} |\geq \frac{1}{\sqrt{n}}\frac{\delta}{4}$, we get

$$\bigg(\int_{\mathbb{R}^n}  |z_{k_a}+ \lambda(a-y)_{k_a} |^\beta\frac{ z_{k_a}}{(1 +|z|^2)^{n+1}}  dz\bigg)^2\geq c_\delta>0.$$Thus,

\begin{eqnarray}
\nonumber  \big<\partial J(u),  \sum_{k=1}^n b_k \int_{\mathbb{R}^n}  |z_k+ \lambda(a-y)_k |^\beta\frac{ z_k}{(1 +|z|^2)^{n+1}}  dz \; \dis\alpha\frac{1}{\lambda}\frac{\partial \delta_{(a, \lambda)}}{\partial a_k}\big> &\leq& - \frac{ c_1}{\lambda^\beta},\\
\nonumber &\leq& - c \; \Big(\frac{ 1}{\lambda^\beta} + \frac{|\nabla K(a)|}{\lambda}\Big).
\end{eqnarray}

\n Lastly, if $\lambda|a-y|\geq \frac{1}{2\delta}$, by the first expansion of Proposition \ref{lem2}, we have

$$ <\partial J(u), \sum_{k=1}^nb_k \mbox{ sign } (a-y)_k \dis\alpha\frac{1}{\lambda}\frac{\partial \delta_{(a, \lambda)}}{\partial a_k}> \leq - c \;\sum_{k=1}^n b_k^2 \frac{|a-y|^{\beta-1}}{\lambda}$$$$ + O\bigg(\sum_{j=2}^{\tilde{\beta}}\frac{ |a -y|^{\beta-j} }{\lambda^j}\bigg)+ O\bigg(\frac{1}{\lambda^{\beta}}\bigg).$$Oberve that for every $j= 2, \ldots, \tilde{\beta}$

$$\frac{|a-y|^{\beta-j}}{\lambda^j} = o\Big( \frac{|a-y|^{\beta-1}}{\lambda} \Big), \mbox{ as } \delta \mbox{ small enough}.$$
Also,

\begin{equation}\label{3c}
\frac{1}{\lambda^\beta} = o\Big( \frac{|a-y|^{\beta-1}}{\lambda} \Big), \mbox{ as } \delta \mbox{ small enough}.
\end{equation}
Then, we obtain

$$ <\partial J(u), \sum_{k=1}^nb_k \mbox{ sign } (a-y)_k \dis\alpha\frac{1}{\lambda}\frac{\partial \delta_{(a, \lambda)}}{\partial a_k}> \leq - c \frac{ |a -y|^{\beta-1} }{\lambda},$$since $\dis \sum_{k=1}^n|(a -y)_k|^{\beta-1}\sim |a -y|^{\beta-1}$. Now by \eqref{3a1} and \eqref{3c}, we derive from the above inequality that

$$ <\partial J(u), \sum_{k=1}^nb_k \mbox{ sign } (a-y)_k \dis\alpha\frac{1}{\lambda}\frac{\partial \delta_{(a, \lambda)}}{\partial a_k}> \leq - c \Big(\frac{1}{\lambda^{\beta}}+ \frac{|\nabla K (a)|}{\lambda}\Big).$$
Hence claim \eqref{aa} follows.

\smallskip

\n\textbf{$\bullet$ Pseudo-gradient in $\widetilde{V}_2(1, \varepsilon)$:}

\n Let $W_2$ be the following vector field. $\forall u= \alpha\delta_{(a, \lambda)}\in \widetilde{V}_2(1, \varepsilon)$,
$$W_2(u)=  \theta_1(\lambda |a-y|)\dis \alpha\lambda \frac{ \partial \delta_{(a, \lambda)}}{\partial \lambda}+ \theta_2(\lambda |a-y|) \sum_{k=1}^n b_k \dis \int_{\mathbb{R}^n}  \frac{ z_{k} |z_{k}+ \lambda(a-y)_{k} |^\beta}{(1 +|z|^2)^{n+1}}  dz$$$$+ \theta_3(\lambda |a-y|) \sum_{k=1}^n b_k \mbox{ sign } (a-y)_k \; \dis\alpha\frac{1}{\lambda}\frac{\partial \delta_{(a, \lambda)}}{\partial a_k}.$$
Observe that, if $\lambda |a-y| \leq {\delta}$, by  the expansion  of Proposition \ref{lem1}, we get

$$
<\partial J(u),   \dis \alpha\lambda \frac{ \partial \delta_{(a, \lambda)}}{\partial \lambda}>\leq  - c \; \frac{ H(a, a)}{\lambda^{n-4}} + O(|a-y|^{n-4})
$$$$\leq -\frac{c}{\lambda^{n-4}},$$
since $|a-y|^{n-4}= o\Big(\frac{1}{\lambda^{n-4}}\Big)$ as $\delta$ small enough and $H$ is a positive regular function on $\Omega^2$. Thus by \eqref{3b}, we obtain

$$
<\partial J(u),   \dis \alpha\lambda \frac{ \partial \delta_{(a, \lambda)}}{\partial \lambda}>\leq  - c \Big(\frac{1}{\lambda^{n-4}}+ \frac{|\nabla K(a)|}{\lambda}\Big).$$
If $\lambda |a-y| \geq \frac{\delta}{2}$, we proceed exactly as in $\widetilde{V}_1(1, \varepsilon)$. We therefore obtain

\begin{equation}\label{3.7}
 <\partial J(u), W_2(u)> \leq - c \Big(\frac{1}{\lambda^{n-4}}+ \frac{|\nabla K(a)|}{\lambda}\Big).
\end{equation}

\smallskip

\n\textbf{$\bullet$ Pseudo-gradient in $\widetilde{V}_3(1, \varepsilon)$:}

\n Let $W_3$ be the following vector field. $\forall u= \alpha\delta_{(a, \lambda)}\in \widetilde{V}_3(1, \varepsilon)$,
$$W_3(u)=  \theta_1(\lambda^{n-4} |a-y|^\beta)\dis \alpha\lambda \frac{ \partial \delta_{(a, \lambda)}}{\partial \lambda}+ (1-\theta_1)(\lambda^{n-4} |a-y|^\beta) \sum_{k=1}^n b_k \dis  \mbox{ sign } (a-y)_k \; \dis\alpha\frac{1}{\lambda}\frac{\partial \delta_{(a, \lambda)}}{\partial a_k}.$$
We claim that

\begin{equation}\label{3k}
 <\partial J(u), W_3(u)> \leq - c \Big(\frac{1}{\lambda^{\min(n-1, \beta)}}+ \frac{|\nabla K(a)|}{\lambda}\Big).
 \end{equation}
Indeed,  if $\lambda^{n-4} |a-y|^\beta \leq {\delta}$ in   the expansion  of Proposition \ref{lem1}, we have

$$O\Big(|a-y|^\beta\Big) = o\Big(\frac{1}{\lambda^{n-4}}\Big),  \qquad \mbox{ taking } \delta \mbox{ small enough, } $$
Thus,

$$
 <\partial J(u), \dis \alpha\lambda \frac{ \partial \delta_{(a, \lambda)}}{\partial \lambda}> \leq - \frac{c}{\lambda^{n-4}}.
$$Moreover,

$$\frac{|\nabla K(a)| }{\lambda} = O\Big(\frac{|a-y|^{\beta -1}}{{\lambda}}\Big)= O\Big(\frac{1}{\lambda^{1+ \frac{(n-4)(\beta-1)}{\beta}}}\Big)= o\Big(\frac{1}{\lambda^{n-4}}\Big)$$Therefore, we can appear $\dis \frac{\nabla K(a)}{\lambda}$ in the latest upper bound. Hence

$$
 <\partial J(u), \dis \alpha\lambda \frac{ \partial \delta_{(a, \lambda)}}{\partial \lambda}> \leq - c\bigg( \frac{1}{\lambda^{n-2}}+ \frac{|\nabla K(a)|}{\lambda}\bigg)
$$$$\leq - c\bigg( \frac{1}{\lambda^{\min(n-1, \beta)}}+ \frac{|\nabla K(a)|}{\lambda}\bigg).$$
Now, if $\lambda^{n-4} |a-y|^\beta >\frac{\delta}{2}$, by the first expansion of Proposition \ref{lem2}, we have

$$ <\partial J(u), \sum_{k=1}^nb_k \mbox{ sign } (a-y)_k \dis\alpha\frac{1}{\lambda}\frac{\partial \delta_{(a, \lambda)}}{\partial a_k}> \leq - c \;\sum_{k=1}^n b_k^2 \frac{|a-y|^{\beta-1}}{\lambda} + O\bigg(\sum_{j=2}^{\tilde{\beta}}\frac{ |a -y|^{\beta-j} }{\lambda^j}\bigg)$$$$+ O\bigg(\frac{1}{\lambda^{\min(\beta, n)}}\bigg)+ O\Big(\frac{1}{\lambda^{n-1}}\Big).$$Oberve that for every $j= 2, \ldots, n$

$$\frac{|a-y|^{\beta-j}}{\lambda^j} = o\Big( \frac{|a-y|^{\beta-1}}{\lambda} \Big), \mbox{ as } \lambda \mbox{ goes to } +\infty.$$
Indeed,

$$\frac{|a-y|^{\beta-j}}{\lambda^j} \frac{\lambda}{|a-y|^{\beta-1}}= \frac{1}{(\lambda|a-y|)^{j-1}}\leq (\frac{2}{\delta})^\frac{j-1}{\beta}\frac{1}{\lambda^{(j-1)(1-\frac{n-4}{\beta})}},$$which goes to zero when $\lambda$ goes to $+\infty$. In addition,

\begin{equation}\label{3e}
\frac{1}{\lambda^{n-1}} = o\Big( \frac{|a-y|^{\beta-1}}{\lambda} \Big), \mbox{ and } \frac{1}{\lambda^{\min(\beta, n)}}=  o\Big( \frac{|a-y|^{\beta-1}}{\lambda} \Big) \mbox{ as } \lambda\rightarrow +\infty.
\end{equation}
Therefore,

$$ <\partial J(u), \sum_{k=1}^nb_k \mbox{ sign } (a-y)_k \dis\alpha\frac{1}{\lambda}\frac{\partial \delta_{(a, \lambda)}}{\partial a_k}> \leq - c  \frac{|a-y|^{\beta-1}}{\lambda}$$and by \eqref{3a} and \eqref{3e}, we obtain

\begin{equation}\label{}
 <\partial J(u), \sum_{k=1}^nb_k \mbox{ sign } (a-y)_k \dis\alpha\frac{1}{\lambda}\frac{\partial \delta_{(a, \lambda)}}{\partial a_k}> \leq - c \Big(\frac{1}{\lambda^{\min(n-1, \beta)}}+ \frac{|\nabla K (a)|}{\lambda}\Big).
 \end{equation}
Hence our claim \eqref{3k} follows.

\smallskip

\begin{pfn}{\bf{Theorem \ref{th3}.}}
Let  $W$ the vector field in $\widetilde{V}(1, \varepsilon)$   defined by convex combination of $W_1, W_2$ and $W_3$. By \eqref{aa}, \eqref{3.7} and \eqref{3k}, we have

$$
 <\partial J(u), \dis W(u)> \leq - c\bigg( \frac{1}{\lambda^{\min(n-1, \beta)}}+ \frac{|\nabla K(a)|}{\lambda}\bigg).
$$

\n In the above construction of $W$, we observe that the Palais-Smale condition is satisfied along  the decreasing flow lines of the pseudo-gradient $W$ as long as the concentration points $a(s)$ of the flow do not enter in some neighborhood of any critical point $\dis y \in \mathcal{K}^+_{< n-4}\cup \mathcal{K}^+_{n-4}\cup \mathcal{K}_{>n-4}$, since $\lambda(s)$ decreases on the flow line in this region. However, if $a(s)$ is near a critical point  $\dis y \in \mathcal{K}^+_{< n-4}\cup \mathcal{K}^+_{n-4}\cup \mathcal{K}_{>n-4}$, $\lambda(s)$ increases on the flow line and goes to $+\infty$.  Thus, we obtain a critical point at infinity.\\
In this statement, the functional $J$ can be expended  after a suitable change of variables as

$$J(\alpha P \delta_{a, \lambda} + \bar{v}) = \dis J(\widetilde{\alpha} P \delta_{\widetilde{a}, \widetilde{\lambda}})= \dis\frac{S_n}{\widetilde{\alpha}^\frac{8}{n-4} (K(\widetilde{a}))^\frac{n-4}{4}}\Big(1+\frac{1}{\widetilde{\lambda}^\beta}\Big).$$
Thus, the index of such critical point at infinity is $n-\widetilde{i}(y)$. Since $J$ behaves in this region as $ \dis\frac{1}{ K^\frac{n-4}{4}}$. This conclude the proof of Theorem \ref{th3}.

\end{pfn}

In the second part of this section, we focus on $\widetilde{V}(p, \varepsilon), p\geq 2$. We characterize the critical points at infinity in these sets in order to give a complete description of the loss of compactness of problem \eqref{1.1} under $(f)_{\beta}$-condition, where $\beta\in (1, n-4]$. Let

\begin{eqnarray}
\nonumber \widetilde{V}_{< n-4}(p, \varepsilon) &=& \{u= \dis\sum_{i=1}^p\alpha_i\delta_{(a_i, \lambda_i)} + v\in \widetilde{V}(p, \varepsilon), \; a_i\in B(y_{\ell_i}, \rho), \;  y_{\ell_i}\in \mathcal{K}\\ & & \nonumber \mbox{ with } \beta= \beta(y_{\ell_i}) < n-4, \forall i=1, \ldots, p\},\\
\nonumber \widetilde{V}_{n-4}(p, \varepsilon) &=& \{u= \dis\sum_{i=1}^p\alpha_i\delta_{(a_i, \lambda_i)} + v\in \widetilde{V}(p, \varepsilon), \; a_i\in \nonumber B(y_{\ell_i}, \rho), \;  y_{\ell_i}\in \mathcal{K}\\ && \nonumber \mbox{ with } \beta= \beta(y_{\ell_i}) = n-4, \forall i=1, \ldots, p\}.
\end{eqnarray}

\n We introduce now the following two Lemmas.

\begin{lem}\label{lem3}
There exists a pseudo-gradient $W_1$ in $\widetilde{V}_{< n-4}(p, \varepsilon)$ such that for any $u= \dis \sum_{i=1}^p\alpha_i P \delta_{a_i, \lambda_i}  \in \widetilde{V}_{< n-4}(p, \varepsilon)$, we have

$$<\partial J(u), W_1(u)> \leq -c \bigg(\sum_{i=1}^p\Big(\frac{1}{\lambda_i^{\beta}}+ \frac{|\nabla K(a_i)|}{\lambda_i}\Big)+ \sum_{j\neq i}\varepsilon_{ij}\bigg).$$
Moreover, the only situation when the  $\lambda_i(s), i=1, \ldots, p, s\geq 0,$ are not bounded is when  $a_i(s)$ goes  to $\dis y_{\ell_i} \in \mathcal{K}^+_{< n-4}, \forall i=1, \ldots, p$ with $y_{\ell_i} \neq y_{\ell_j}, \forall i\neq j$.
\end{lem}

\begin{lem}\label{lem4}
There exists a pseudo-gradient $W_2$ in $\widetilde{V}_{n-4}(p, \varepsilon)$ such that for any $u= \dis \sum_{i=1}^p\alpha_i P \delta_{a_i, \lambda_i}  \in \widetilde{V}_{n-4}(p, \varepsilon)$, we have

$$<\partial J(u), W_2(u)> \leq -c \bigg(\sum_{i=1}^p\Big(\frac{1}{\lambda_i^{n-4}}+ \frac{|\nabla K(a_i)|}{\lambda_i}\Big)+ \sum_{j\neq i}\varepsilon_{ij}\bigg).$$
Moreover, the only situation when the  $\lambda_i(s), i=1, \ldots, p, s\geq 0,$ are not bounded is when  $a_i(s)$ goes  to $\dis y_{\ell_i} \in \mathcal{K}^+_{ n-4}, \forall i=1, \ldots, p$ with $y_{\ell_i} \neq y_{\ell_j}, \forall i\neq j$ and $\rho(y_{\ell_1}, \ldots, y_{\ell_p})>0$.
\end{lem}

The proof of Lemmas  \ref{lem3} and \ref{lem4} will be given at the end of this section. We now state  the proof of Theorem \ref{th3}.

\smallskip

\begin{pfn}{\bf{Theorem \ref{th3}.}}
It follows from the following Lemma.

\begin{lem}\label{lem5} Under the assumption that $K$ satisfies  $(B)$ and $(f)_{\beta}$ for $\beta\in (1, n-4]$, there exists a pseudo-gradient $W$ in $\widetilde{V}(p, \varepsilon)$ such that for any $u= \dis \sum_{i=1}^p\alpha_i P \delta_{a_i, \lambda_i}  \in \widetilde{V}(p, \varepsilon)$, we have

\begin{equation}\label{3.10}
<\partial J(u), \widetilde{W}(u)> \leq -c \bigg(\sum_{i=1}^p\Big(\frac{1}{\lambda_i^{\beta}}+ \frac{|\nabla K(a_i)|}{\lambda_i}\Big)+ \sum_{j\neq i}\varepsilon_{ij}\bigg).\end{equation}
Moreover, the only case where $\lambda_i(s), i=1, \ldots, p, s\geq 0$ are not bounded is when $a_i(s)$ goes to $y_{\ell_i}\in  \mathcal{K}^+_{< n-4}\cup \mathcal{K}^+_{n-4}, \forall i=1, \ldots, p$ with $y_{\ell_i}\neq y_{\ell_j}, \forall i\neq j$ and $(y_{\ell_1}, \ldots, y_{\ell_p})\in  \mathcal C_{n-4}^{\infty}\cup \mathcal C_{<n-4}^{\infty}\cup \Big(\mathcal C_{n-4}^{\infty}\times\mathcal C_{<n-4}^{\infty}\Big)$.

\end{lem}
\smallskip

\begin{pfn}{\bf{Lemma \ref{lem5}.}}
Let $u= \dis \sum_{i=1}^p\alpha_i P \delta_{a_i, \lambda_i}  \in \widetilde{V}(p, \varepsilon), p\geq 2$. By Lemmas \ref{lem3} and \ref{lem4}, it remains only to consider the case where  $$u= \dis \sum_{i \in I_1}\alpha_i P \delta_{a_i, \lambda_i}+ \dis \sum_{i \in I_2}\alpha_i P \delta_{a_i, \lambda_i},$$
with $I_1\neq\emptyset$, $I_2\neq\emptyset$, $u_1:=  \dis \sum_{i \in I_1}\alpha_i P \delta_{a_i, \lambda_i}\in  \widetilde{V}_{<n-4}( \sharp I_1, \varepsilon)$ and  $u_2:=  \dis \sum_{i \in I_2}\alpha_i P \delta_{a_i, \lambda_i}\in  \widetilde{V}_{n-4}( \sharp I_2, \varepsilon)$. We order the $\lambda_i$'s, we can assume that $$\lambda_1\leq \lambda_2\leq \cdots \leq\lambda_p.$$
Three cases may occur.

\n$\bullet$ First case: $u_1\not\in \Big\{u= \dis \sum_{j=1}^{\sharp I_1}\alpha_j P \delta_{a_j, \lambda_j} \in  \widetilde{V}_{<n-4}( \sharp I_1, \varepsilon)$ with $\dis y_{\ell_j} \in \mathcal{K}^+_{<n-4}$, $y_{\ell_j} \neq y_{\ell_k}, \forall 1\leq j\neq k\leq \sharp I_1\Big\}$.

\n Let $\widetilde{W}_1(u)= W_1(u_1)$ where $W_1$ is the pseudo-gradient defined in Lemma \ref{lem3}. Observe that the maximum of the $\lambda_i(s), i\in I_1$ does not increase through $W_1$. Moreover, by Lemma \ref{lem3}, we have

\begin{equation}\label{c1c1}
<\partial J(u), \widetilde{W}_1(u)> \leq -c \bigg(\sum_{i \in I_1}\Big(\frac{1}{\lambda_i^{\beta}}+ \frac{|\nabla K(a_i)|}{\lambda_i}\Big)+ \sum_{i, j \in I_1, j\neq i}\varepsilon_{ij}\bigg)\end{equation}
$$+ O\bigg(\sum_{i\in I_1, j \not \in I_1}\varepsilon_{ij}\bigg).$$
Let $k_0$ be an index in $I_1$ such that

$$\lambda_{k_0}^{\beta(y_{\ell_{k_0}})}= \min \{\lambda_i^{\beta(y_{\ell_i})}, i\in I_1\}.$$
Define $$\widetilde{I}_1 = \{i, 1\leq i \leq p, s.t. \lambda_i^{\beta(y_{\ell_i})}\geq \dis \frac{1}{10}\lambda_{k_0}^{\beta(y_{\ell_{k_0}})}\}.$$
Observe that $I_1 \subset \widetilde{I}_1$. Our first goal is to make appears in the upper bound of \eqref{c1c1} all indices $i\in \widetilde{I}_1$. For each index $i$ we define the following vector field.

\begin{equation}\label{x}
    X_i(u)=\sum_{k=1}^n b_k \mbox{ sign } (a_i-y_{\ell_i})_k \frac{1}{\lambda_{i}} \dis \frac{\partial
P\delta_{a_i, \lambda_i}}{\partial (a_{i})_k}.
\end{equation}
Using the first expansion of Lemma \ref{lem2}, we have

\begin{equation}\label{c2c2}
<\partial J(u), X_i(u)> \leq  - \frac{c_3}{K(a_i)} \alpha_i^2J(u) \sum_{k=1}^n b_k^2   \frac{|a_i-y_{\ell_i}|^{\beta-1}}{\lambda_i}
 \end{equation}$$+     O\Big(\sum_{j=2}^{[\beta]}\frac{|a_i- y_{\ell_i}|^{\beta-j}}{\lambda_i^j}\Big) +    O\Big(\frac{1}{\lambda_i^{\beta}}\Big) + o\Big( \dis \sum_{j \neq i} \varepsilon_{ij}   \Big),$$
since $\dis \Big|\frac{1}{\lambda_i}\frac{\partial \varepsilon_{i j}}{\partial a_i}\Big|= o(\varepsilon_{ij})$ if $|a_i-a_j|>\rho$. Let $M$ be a large positive constant. Observe that if $\dis| \lambda_i(a_i- y_{\ell_i}) |\leq M$,

$$\frac{|a_i- y_{\ell_i}|^{\beta-j}}{\lambda_i^j}= O\Big(\frac{1}{\lambda_i^\beta}\Big),$$and if $\dis| \lambda_i(a_i- y_{\ell_i}) |\geq M$,

$$\frac{|a_i- y_{\ell_i}|^{\beta-j}}{\lambda_i^j}= o\Big(\frac{|a_i- y_{\ell_i}|^{\beta-1}}{\lambda_i}\Big),$$ and $$\frac{1}{\lambda_i^\beta}= o\Big(\frac{|a_i- y_{\ell_i}|^{\beta-1}}{\lambda_i}\Big) \mbox{ as M large}.$$
Therefore, for $m_1>0$ very small, we get from \eqref{c1c1} and \eqref{c2c2}

\begin{equation}\label{}
<\partial J(u), \widetilde{W}_1(u) + m_1\sum_{i\in \widetilde{I}_1\setminus I_1}X_i(u)> \leq -c \bigg(\sum_{i \in \widetilde{I}_1}\Big(\frac{1}{\lambda_i^{\beta}}+ \frac{|\nabla K(a_i)|}{\lambda_i}\Big)+ \sum_{i, j \in I_1, j\neq i}\varepsilon_{ij}\bigg)\end{equation}
$$+ o\Big(\sum_{j=1}^p\frac{1}{\lambda_j^\beta}\Big),$$since under $(f)_\beta$-condition $|\nabla K(a_i)|= O(|a_i- y_{\ell_i}|^{\beta-1})$ and $$\varepsilon_{ij}= O\Big(\frac{1}{(\lambda_i\lambda_j)^\frac{n-2}{2}}\Big)= o\Big(\frac{1}{\lambda_i^\beta}+ \frac{1}{\lambda_j^\beta}\Big), \forall i\in I_1, j\not\in I_1.$$

\n In order to appear $\dis -\sum_{i, j \in \widetilde{I}_1, j\neq i}\varepsilon_{ij}$ we will decrease all the $\lambda_i, i\in \widetilde{I}_1\setminus I_1$ with different speed. Let $Z_1(u)= \dis \sum_{i\in \widetilde{I}_1\setminus I_1}-2^i \lambda_{i} \dis \frac{\partial
P \delta_{a_i, \lambda_i}}{\partial \lambda_{i}}$. Observe that

$$2^i\lambda_i \frac{\partial  \varepsilon_{ij}}{\partial\lambda_i} + 2^j\lambda_j \frac{\partial  \varepsilon_{ij}}{\partial\lambda_j}\leq -c\varepsilon_{ij}, \forall i<j.$$
Thus, using the first expansion of Proposition \ref{lem1}, we get

$$\Big \langle \partial J(u), Z_1(u)\Big\rangle \leq
 - c \dis\sum_{i, j \in \widetilde{I}_1\setminus I_1, j\neq i}  \varepsilon_{i j}+ O\Big( \sum_{i\in \widetilde{I}_1\setminus I_1}\frac{1}{\lambda_{i}^{\beta}} + \sum_{j=2}^{[\beta]}\frac{|a_i- y_{\ell_i}|^{\beta-j}}{\lambda_i^j}\Big).$$
Therefore, for $m_2>2$ very small, we obtain

\begin{equation}\label{}
<\partial J(u), \widetilde{W}_1(u) +m_2Z_1(u) + m_1\sum_{i\in \widetilde{I}_1\setminus I_1}X_i(u)> \leq -c \bigg(\sum_{i \in \widetilde{I}_1}\Big(\frac{1}{\lambda_i^{\beta}}+ \frac{|\nabla K(a_i)|}{\lambda_i}\Big)+ \sum_{i, j \in \widetilde{I}_1, j\neq i}\varepsilon_{ij}\bigg)\end{equation}
$$+ o\Big(\sum_{j=1}^p\frac{1}{\lambda_j^\beta}\Big).$$
Now let $R$ be the set of the remainder indices.  So $R= \{i\in I_2, \lambda_i^\beta < \frac{1}{10}\lambda_{k_0}^{\beta}\}$ and denote $\widetilde{u}= \dis \sum_{i\in R}\alpha_i P \delta_{a_i, \lambda_i}$. Observe that $\widetilde{u}\in \widetilde{V}_{n-2}(\sharp R, \varepsilon)$, therefore, we can apply the associated vector field $W_2(\widetilde{u})$ defined in Lemma \ref{lem4}. For $\widetilde{W}_2(u)= W_2(\widetilde{u})$ we get by Lemma \ref{lem4}

$$<\partial J(u), \widetilde{W}_2(u)> \leq -c \bigg(\sum_{i\in R}\Big(\frac{1}{\lambda_i^{\beta}}+ \frac{|\nabla K(a_i)|}{\lambda_i}\Big)+ \sum_{i,j \in R, j\neq i}\varepsilon_{ij}\bigg)$$$$+ O\bigg(\sum_{i\in R, j \not \in R}\varepsilon_{ij}\bigg).$$
We let in this case $W= \widetilde{W}_1+\widetilde{W}_2+ m_2Z_1+ m_1 \sum_{i\in \widetilde{I}_1\setminus I_1}X_i(u)$. It satisfies

\begin{equation}\label{}
<\partial J(u), W(u)> \leq -c \bigg(\sum_{i=1}^p\Big(\frac{1}{\lambda_i^{\beta}}+ \frac{|\nabla K(a_i)|}{\lambda_i}\Big)+ \sum_{ j\neq i}\varepsilon_{ij}\bigg).
\end{equation}

\n $\bullet$ Second case:  $u_2\not\in \Big\{u= \dis \sum_{j=1}^{\sharp I_2}\alpha_j P \delta_{a_j, \lambda_j} \in  \widetilde{V}_{n-4}( \sharp I_2, \varepsilon)$ with $\dis y_{\ell_j} \in \mathcal{K}^+_{n-4}$, $y_{\ell_j} \neq y_{\ell_k}, \forall 1\leq j\neq k\leq \sharp I_2$ and $\rho(y_{\ell_1}, \ldots, y_{\ell_p})>0 \Big\}$.

\n Let  $\widetilde{W}_2(u)= W_2(u_2)$ where $W_2$ is defined in Lemma \ref{lem4}. We then have:

$$<\partial J(u), \widetilde{W}_2(u)> \leq -c \bigg(\sum_{i\in I_2}\Big(\frac{1}{\lambda_i^{\beta}}+ \frac{|\nabla K(a_i)|}{\lambda_i}\Big)+ \sum_{i,j \in I_2, j\neq i}\varepsilon_{ij}\bigg)$$$$+ O\bigg(\sum_{i\in I_2, j \not \in I_2}\varepsilon_{ij}\bigg).$$
As in the first case, we denote by $k_0$ the index of $I_2$ satisfying
$$\lambda_{k_0}= \min \{\lambda_i, i\in I_2\}$$
and we define $$\widetilde{I}_2 = \{i, 1\leq i \leq p, s.t. \lambda_i^{\beta(y_{\ell_i})}\geq \dis \frac{1}{10}\lambda_{k_0}^{n-4}\}$$
and $E= \widetilde{I}_2^c$. Let $$Z_2(u)= \dis \sum_{i\in \widetilde{I}_2\setminus I_2}-2^i \lambda_{i} \dis \frac{\partial
P \delta_{a_i, \lambda_i}}{\partial \lambda_{i}},$$$$\widetilde{W}_1(u)= W_1\Big(\dis\sum_{i\in E}\alpha_iP \delta_{a_i, \lambda_i}\Big),$$and $X_i(u)$ the vector field defined in \eqref{x}. By the same computation of the first section, we get for $W= \widetilde{W}_2 + \widetilde{W}_1+ m_2Z_2+ m_1\sum_{i\in \widetilde{I}_2\setminus I_2}X_i(u)$,

$$
<\partial J(u), W(u)> \leq -c \bigg(\sum_{i=1}^p\Big(\frac{1}{\lambda_i^{\beta}}+ \frac{|\nabla K(a_i)|}{\lambda_i}\Big)+ \sum_{ j\neq i}\varepsilon_{ij}\bigg).
$$

\n $\bullet$ Third case:  $u_1 \in \Big\{u= \dis \sum_{j=1}^{\sharp I_1}\alpha_j P \delta_{a_j, \lambda_j} \in  \widetilde{V}_{<n-4}( \sharp I_1, \varepsilon)$ with $\dis y_{\ell_j} \in \mathcal{K}^+_{<n-4}$ and $y_{\ell_j} \neq y_{\ell_k}, \forall 1\leq j\neq k\leq \sharp I_2\Big\}$ and $u_2\in \Big\{u= \dis \sum_{j=1}^{\sharp I_2}\alpha_j P \delta_{a_j, \lambda_j} \in  \widetilde{V}_{n-4}( \sharp I_2, \varepsilon)$ with $\dis y_{\ell_j} \in \mathcal{K}^+_{n-4}$, $y_{\ell_j} \neq y_{\ell_k}, \forall 1\leq j\neq k\leq \sharp I_2$ and $\rho(y_{\ell_1}, \ldots, y_{\ell_p})>0 \Big\}$.

\n Let $\widetilde{W}_1(u)= W_1(u_1)$ and $\widetilde{W}_2(u)= W_2(u_2)$ where $W_1$ and $W_2$ are the vector fields defined in Lemmas \ref{lem3} and \ref{lem4} respectively. Using the above estimates, we get for $W= \widetilde{W}_1+ \widetilde{W}_2$

$$
<\partial J(u), W(u)> \leq -c \bigg(\sum_{i=1}^p\Big(\frac{1}{\lambda_i^{\beta}}+ \frac{|\nabla K(a_i)|}{\lambda_i}\Big)+ \sum_{ j\neq i}\varepsilon_{ij}\bigg).
$$
This finishes the proof of Lemma \ref{lem5} and then the proof of Theorem \ref{th3} follows.
\end{pfn}

\end{pfn}

\begin{pfn}{\bf{Lemma \ref{lem3}.}}
We divide $\widetilde{V}_{1}(p, \varepsilon)$ as follows. Let $\delta>0$ and small.

\n $W_{1}(p,\varepsilon) := \Big\{u=\dis \sum_{j=1}^{p}\alpha_{j}
P \delta_{a_{j},  \lambda_{j}}\in \widetilde{V}_{1}(p,\varepsilon),\;
y_{l_{j}}\neq y_{l_{k}} ,\; \forall j \neq k,\;-
\sum_{k=1}^{n}b_{k}(y_{l_{j}}) > 0, \mbox{and
}\lambda_{j}|a_{j}-y_{l_{j}}|< \delta,\; \forall j=1,...,p \Big\}.$\\
$W_{2}(p,\varepsilon) := \Big\{u=\dis \sum_{j=1}^{p}\alpha_{j}
P \delta_{a_{j}, \lambda_{j})}\in \widetilde{V}_{1}(p,\varepsilon),\;
y_{l_{j}}\neq y_{l_{k}}, \; \forall j \neq
k,\;\lambda_{j}|a_{j}-y_{l_{j}}|< \delta,\; \forall j=1,..,p
\mbox{ and there exist at least} j_{1} \mbox{ such
that }-  \sum_{k=1}^{n}b_{k}(y_{l_{j_{1}}}) < 0 \Big\}.$\\
$
W_{3}(p,\varepsilon) := \Big\{u=\dis \sum_{j=1}^{p}\alpha_{j}
P \delta_{a_{j}, \lambda_{j}}\in \widetilde{V}_{1}(p,\varepsilon),\;
y_{l_{j}}\neq y_{l_{k}}, \; \forall j \neq k,\; \mbox{and there
exist at least } j_1, s, t,\newline  \lambda_{j_1}|a_{j_1}-y_{l_{j_1}}|\geq\dis \frac{\delta}{2}  \Big\}.
$\\
$
W_{4}(p,\varepsilon) := \Big\{u=\dis \sum_{j=1}^{p}\alpha_{j}
P \delta_{a_{j}, \lambda_{j}}\in \widetilde{V}_{1}(p,\varepsilon),\;
 \mbox{such that there exist } j\neq k
  \mbox{ with }y_{l_{j}}= y_{l_{k}} \Big\}.
$\\

\n $\bullet$ Pseudo-gradient in $W_1(p, \varepsilon)$. In this region, we have $ |a_{j}-y_{l_{j}}|^\beta= o \dis\Big(\frac{1}{\lambda_{j}^\beta}\Big), \forall j=1, \ldots, p$ and $\dis \frac{1}{(\lambda_j\lambda_k)^\frac{n-2}{2}} = o\Big(\frac{1}{\lambda_{j}^\beta}\Big) + o\Big(\frac{1}{\lambda_{k}^\beta}\Big)$ since $\beta< n-2$. Thus, for $V_1(u)= \dis \sum_{j=1}^{p}  \lambda_{j} \dis \frac{\partial
P \delta_{a_j, \lambda_j}}{\partial \lambda_{j}}$, we have:

$$
<\partial J(u), V_1(u)> \leq -c \bigg(\sum_{i=1}^p\Big(\frac{1}{\lambda_i^{\beta}}+ \frac{|\nabla K(a_i)|}{\lambda_i}\Big)+ \sum_{ j\neq i}\varepsilon_{ij}\bigg),
$$since $\dis  \frac{|\nabla K(a_i)|}{\lambda_i}= o \dis\Big(\frac{1}{\lambda_{j}^\beta}\Big)$.

\n $\bullet$ Pseudo-gradient in $W_2(p, \varepsilon)$. Let $k_0$ an index such that $\lambda_{k_0}^\beta= \min \{\lambda_j^\beta,$ such that $-\sum_{k=1}^{n}b_{k}(y_{l_{j}}) <0\}$. We denote by $I_1=\{j=1, \ldots, p, $ such that, $\lambda_j^\beta\leq \dis \frac{1}{2} \lambda_{k_0}^\beta\}$. Observe that $\widetilde{u} :=   \sum_{j\in I_1}\alpha_{j} P \delta_{a_{j}, \lambda_{j}}\in W_1(\sharp I_1, \varepsilon)$. Using the same previous technics, we get for $\widetilde{V}_1(u)= V_1(\widetilde{u})$,

$$
<\partial J(u), \widetilde{V}_1(u)- \sum_{j\not\in I_1}\lambda_{j} \dis \frac{\partial
P \delta_{a_j, \lambda_j}}{\partial \lambda_{j}}> \leq -c \bigg(\sum_{i=1}^p\Big(\frac{1}{\lambda_i^{\beta}}+ \frac{|\nabla K(a_i)|}{\lambda_i}\Big)+
\sum_{ j\neq i}\varepsilon_{ij}\bigg).
$$

\n $\bullet$ Pseudo-gradient in $W_3(p, \varepsilon)$. Let $k_1$ an index such that $\lambda_{k_1}^\beta= \min \{\lambda_j^\beta,$ such that $\lambda_{j}|a_{j}-y_{l_{j}}|\geq\dis \frac{\delta}{2}\}$ and let $J_1= \{j=1, \ldots, p$ such that $\lambda_j^\beta \leq  \frac{1}{2}\lambda_{k_1}^\beta\}$. For any $j\in J_1$, we have $\lambda_{j}|a_{j}-y_{l_{j}}|\leq\dis \frac{\delta}{2}$. Let $\widetilde{u}:= \sum_{j\in J_1}\alpha_{j} P \delta_{a_{j}, \lambda_{j}}$, $\widetilde{u}\in \dis W_1(\sharp J_1, \varepsilon)\cup W_2(\sharp J_1, \varepsilon)$. We denote by $\widetilde{W}_3(u)=V( \widetilde{u})$ where $V$ is the associated vector field to the above two regions. Using the second expansion of Lemma \ref{lem2} and the previous technics, we get
$$
<\partial J(u), \widetilde{W}_3(u)+ \sum_{i\not\in J_1}\Big(X_i(u)- \lambda_{i} \dis \frac{\partial
P \delta_{a_i, \lambda_i}}{\partial \lambda_{i}}\Big)> \leq -c \bigg(\sum_{i=1}^p\Big(\frac{1}{\lambda_i^{\beta}}+ \frac{|\nabla K(a_i)|}{\lambda_i}\Big)+ \sum_{ j\neq i}\varepsilon_{ij}\bigg).
$$
Here $X_i(u)= \dis\alpha_i\sum_{k=1}^n\frac{b_k}{\lambda_i}\int_{\mathbb{R}^n}  \frac{|x_k+ \lambda_i(a_i-y_{\ell_i})_k |^\beta}{(1+ \lambda_i|(a_i-y_{\ell_i})_k |^{\beta-1}}\frac{ x_k}{(1 +|x|^2)^{n+1}}  dx$ if $\lambda_i|(a_i-y_{\ell_i})_k |\in [\dis \frac{\delta}{2}, \frac{2}{\delta}]$ and $X_i(u)$ is defined by \eqref{x} otherwise.

\n $\bullet$ Pseudo-gradient in $W_4(p, \varepsilon)$.  For any critical point $y_\ell$ of $K$, we denote $\dis B_\rho= \{j=1, \ldots, p$ such that $a_j\in B(y_\ell, \rho)\}$. In this region, there exists at least $\ell$ such that $\sharp B_\ell\geq 2$. Let

$$J_1= \{\ell, 1\leq \ell \leq\sharp \mathcal{K}, \mbox{ such that } \sharp B_\ell \geq 2\}.$$
For any $\ell\in J_1$, we decrease all $\lambda_j$'s, $j\in B_\ell$ as follows. Let

\begin{eqnarray}
\nonumber \psi: \mathbb{R}&\longrightarrow& \mathbb{R}\\
\nonumber t&\longmapsto& \left\{
                           \begin{array}{ll}
                             1 & \hbox{ if } |t|\geq \dis1\\
                             0 & \hbox{ if } |t| \leq \dis \gamma.
                           \end{array}
                         \right.
\end{eqnarray}
where $\gamma>0$ very small. Define $\overline{\psi}(\lambda_j)= \dis \sum_{i\neq j}\psi\Big(\frac{\lambda_j}{\lambda_i}\Big)$ for $j\in B_k$ and $k\in J_1$. By the first expansion of Lemma \ref{lem1}, we get

$$
<\partial J(u), -\dis\sum_{\ell\in J_1} \sum_{j\in B_\ell} \overline{\psi}(\lambda_{j}) \lambda_{j} \dis \frac{\partial
P \delta_{a_j, \lambda_j}}{\partial \lambda_{j}} > \leq -c \sum_{\ell\in J_1}\sum_{ i\neq j \in B_\ell}\varepsilon_{ij} + + O\bigg(\sum_{i=1}^p\frac{1}{\lambda_i^{\beta}}\bigg).
$$
To obtain the required upper bound, we set

$$I_1:= \{j; \lambda_j|a_i-y_{\ell_j} |\geq \delta\}.$$
If $I_1\neq \emptyset$, we use the above vector field (defined in $W_3(p, \varepsilon)$) and using the expansions of Lemma \ref{lem2}, we obtain

$$
<\partial J(u), -\dis\sum_{\ell\in J_1} \sum_{j\in B_\ell} \overline{\psi}(\lambda_{j}) \lambda_{j} \dis \frac{\partial
P \delta_{a_j, \lambda_j}}{\partial \lambda_{j}} + \sum_{i\in I_1}X_i(u)\Big> \leq -c  \bigg(\sum_{i=1}^p\Big(\frac{1}{\lambda_i^{\beta}}+ \frac{|\nabla K(a_i)|}{\lambda_i}\Big)+ \sum_{ j\neq i}\varepsilon_{ij}\bigg).
$$
If $I_1= \emptyset$, we denote by $I_2$ the set of indices constructed by $1$ and all $j$ such that $\lambda_j\sim \lambda_1$, (of the same order). We write $u=  \sum_{i\in I_2}\alpha_{i} P \delta_{a_{i}, \lambda_{i}}+  \sum_{i \not\in I_2}\alpha_{i} P \delta_{a_{i}, \lambda_{i}}$. Observe that $\widetilde{u}:= \sum_{i\in I_2}\alpha_{i} P \delta_{a_{i}, \lambda_{i}} \in W_k(\sharp I_2, \varepsilon), k=1, 2$ or $3$. We then apply the associated vector field denoted $\widetilde{W}_4(u)$. We obtain

$$<\partial J(u), \widetilde{W}_4(u)> \leq -c \bigg(\sum_{i\in I_2}\Big(\frac{1}{\lambda_i^{\beta}}+ \frac{|\nabla K(a_i)|}{\lambda_i}\Big)+ \sum_{i,j \in I_2, j\neq i}\varepsilon_{ij}\bigg)$$$$+ O\bigg(\sum_{i\in I_2, j \not \in I_2}\varepsilon_{ij}\bigg)$$
and therefore,

$$<\partial J(u), \widetilde{W}_4(u) -\dis\sum_{\ell\in J_1} \sum_{j\in B_\ell} \overline{\psi}(\lambda_{j}) \lambda_{j} \dis \frac{\partial
P \delta_{a_j, \lambda_j}}{\partial \lambda_{j}} > \leq -c \bigg(\sum_{i=1}^p\Big(\frac{1}{\lambda_i^{\beta}}+ \frac{|\nabla K(a_i)|}{\lambda_i}\Big)+ \sum_{j\neq i}\varepsilon_{ij}\bigg).$$
This finishes the proof of Lemma \ref{lem3}.

\end{pfn}

\begin{pfn}{\bf{Lemma \ref{lem4}.}}
The situation here is exactly the one of (\cite{CR2}, Proposition 3.7), so we omit the proof here.

\end{pfn}

\cite{}

\section{Proof of the existence results}

\subsection{Proof of Theorem \ref{th2}}

Using the result of Theorem \ref{th1}, the critical points at infinity of the associated variational problem are in one to one correspondence with the elements $(y_1, \ldots, y_p)_{\infty}$, $(y_1, \ldots, y_p)_{\infty}\in \mathcal{C}^\infty$. For each $(y_1, \ldots, y_p)_{\infty}\in \mathcal{C}^\infty$, we denote by $W_u^\infty (y_1, \ldots, y_p)_{\infty}$; the unstable manifold of the critical points at infinity $(y_1, \ldots, y_p)_{\infty}$. Recall that $i(y_1, \ldots, y_p)_{\infty}$ the index of $(y_1, \ldots, y_p)_{\infty}$ is equal to the dimension of $W_u^\infty (y_1, \ldots, y_p)_{\infty}$. Using now the gradient flow of $(-\partial J)$ to deform $\Sigma^+$. It follows then by deformation Lemma (see \cite{BR}), that

\begin{equation}\label{4.1}
\Sigma^+\simeq \dis \bigcup_{(y_1, \ldots, y_p)_{\infty}\in \mathcal{C}^\infty} W_u^\infty (y_1, \ldots, y_p)_{\infty} \cup \bigcup_{w, \partial J(w)=0} W_u(w),
\end{equation}
where $\simeq$  denotes retracts by deformation. It follows from the above deformation retract that the problem \eqref{1.1} has necessary a solution $w$. Otherwise, it follows from \eqref{4.1} that $$1 = \chi(\Sigma^+) = \sum_{(y_1, \ldots, y_p)_{\infty}\in \mathcal C^\infty }(-1)^{i(y_1, \ldots, y_p)_{\infty}},
$$where $\chi$ denotes the Euler-Poincar´e Characteristic. Such an equality contradicts the assumption of Theorem \ref{th2}

\subsection{Proof of Theorem \ref{th4}}

Let
$$J_1(u)= \dis\frac{1}{\bigg(\dis\int_{\Omega}u^\frac{2n}{n-2} dx\bigg)^\frac{n}{n-2}}, \qquad  u\in \Sigma$$be the Euler Lagrange functional associated to Yamabe problem on $\Omega$. Let $$S=\dis \frac{1}{\bigg(\dis\int_{\mathbb{R}^n}\delta_{a, \lambda}^\frac{2n}{n-2} dx\bigg)^\frac{n}{n-2}}$$be the best Sobolev constant. $S$ does not depend on $a$ and $\lambda$. It is known that

$$S= \inf_{u\in \Sigma} J_1(u)$$ and that the infimum is not achieved, see \cite{BrC}.

\n For $c\in \mathbb{R}$ and for any function $f$ on $\Sigma$, we define

$$f^c=\{u\in \Sigma, \; s.t, f(u)\leq c\}.$$It is easy to see that if $|K-1|_{L^\infty(\Omega)}\leq \varepsilon_0$ for $\varepsilon_0$  small enough, we have

\begin{equation}\label{4.1}
    J^{S+ \frac{S}{4}}\subset J_1^{S+ \frac{S}{2}}\subset J^{S+ \frac{3S}{4}}.
\end{equation}
This is due to the fact that $J(u)= J_1(u) ( 1 +O(|\varepsilon_0|))$.

\n Now let $(y_1, \ldots, y_q)_\infty$ be a critical point at infinity of $q$ masses. It is known that the level of $J$ at $(y_1, \ldots, y_q)_\infty$ is given by $\dis S \bigg(\sum_{k=1}^q\frac{1}{K (y_k)^{(n-2)/2}}\bigg)^{2/n}$, see \cite{bouch1}. Hence goes to $qS$ when $\varepsilon_0$ is close to zero. Therefore, for $|\varepsilon_0|$ small enough, we have: \\

\n All critical points at infinity of $J$ of  $q$-masses, $q\geq 2$ are above $S+ \frac{3}{4}S$, \; \, $(*)$

\n and

\n all critical points at infinity of $J$ of one masse are below  $S+ \frac{S}{4}$. \; $(**)$\\

\n Therefore,

\begin{equation}\label{eq41}
\mbox{the functional } J \mbox{ has no critical points at infinity in } {J}^{(\widetilde{S}_{n}+3 \eta)}_{(\widetilde{S}_{n}+\eta)}.
\end{equation}
\n To prove the existence result, we argue by contradiction and we assume that $J$ has no critical points. It follows from \eqref{eq41} that
\begin{equation}\nonumber
J^{\widetilde{S}_{n}+3 \eta}\simeq
J^{\widetilde{S}_{n}+ \eta},
\end{equation}
where $\simeq$ denotes retracts by deformation. Thus by \eqref{4.1},
we derive that

\begin{equation}\label{eq42}
J_{1}^{\widetilde{S}_{n}+2 \eta} \simeq
J^{\widetilde{S}_{n}+ \eta}.
\end{equation}
Now we use the gradient flow of $(-\partial J)$ to deform $J^{\widetilde{S}_{n}+\eta}$. As mentioned above, the only critical points at infinity of $J$ under the level $\widetilde{S}_{n}+\eta$ are $(y)_{\infty}$, $y\in \mathcal{K}_{< n-4}^+ \; \cup \mathcal{K}_{\geq n-4}$. Thus

\begin{equation}\label{eq44}
 J^{\widetilde{S}_{n}+\eta} \simeq \bigcup_{y\in \mathcal{K}_{< n-4}^+ \; \cup \;\mathcal{K}_{\geq n-4}}W^{\infty}_{u} (y).
 \end{equation}
We apply now the Euler-Poincar\'{e} characteristic of both sides of \eqref{eq44}, we get

$$ \chi( J^{\widetilde{S}_{n}+\eta})=\sum_{y\in \mathcal{K}_{< n-4}^+ \; \cup \; \mathcal{K}_{\geq n-4}}(-1)^{n- \widetilde{i}(y)}.$$
Thus by \eqref{eq42}, we obtain

\begin{equation}\label{4.5}
\chi( J_{1}^{\widetilde{S}_{n}+2\eta})=\sum_{y\in \mathcal{K}_{< n-4}^+ \; \cup \; \mathcal{K}_{\geq n-4}}(-1)^{n- \widetilde{i}(y)}.\end{equation}
It is known that $J_{1}^{\widetilde{S}_{n}+2\eta}$ and $\Omega$ has the same homotopy type. See (\cite{bc2} remark 5). Therefore, from \eqref{4.5} we get

\begin{equation}\label{}
\chi( \Omega)=\sum_{y\in \mathcal{K}_{< n-4}^+ \; \cup \; \mathcal{K}_{\geq n-4}}(-1)^{n- \widetilde{i}(y)}.\end{equation}
Such equality contradicts the
assumption of Theorem \ref{th4}. This complete the proof of Theorem \ref{th4}.

\end{document}